%% file: Regular_policies_in_ADP.tex
\input TEXSHOP_macros_new.tex

\input TEXSHOP_small_baseline.tex
\def\section#1{\goodbreak\vskip 3pc plus 6pt minus 3pt\leftskip=-2pc
   \global\advance\sectnum by 1\eqnumber=1\subsectnum=0%
\global\examplnumber=1\figrnumber=1\propnumber=1\defnumber=1\lemmanumber=1\assumptionnumber=1 \conditionnumber =1%
   \line{\hfuzz=1pc{\hbox to 3pc{\bf 
   \vtop{\hfuzz=1pc\hsize=38pc\hyphenpenalty=10000\noindent\uppercase{\the\sectnum.\quad #1}}\hss}}
			\hfill}
			\leftskip=0pc\nobreak\tenf
			\vskip 1pc plus 4pt minus 2pt\noindent\ignorespaces}
\def\subsection#1{\noindent\leftskip=0pc\tenf
   \goodbreak\vskip 1pc plus 4pt minus 2pt
               \global\advance\subsectnum by 1
   \line{\hfuzz=1pc{\hbox to 3pc{\bf \the\sectnum.\the\subsectnum.
   \vtop{\hfuzz=1pc\hsize=38pc\hyphenpenalty=10000\noindent{\bf #1}}\hss}}
                        \hfill}
   \leftskip=0pc\nobreak\tenf
                        \vskip 1pc plus 4pt minus 2pt\nobreak\noindent\ignorespaces}



\def\texshopbox#1{\boxtext{462pt}{\vskip-1.5pc\pshade{\vskip-1.0pc#1\vskip-2.0pc}}}
\def\texshopboxnt#1{\boxtextnt{462pt}{\vskip-1.5pc\pshade{\vskip-1.5pc#1\vskip-2.0pc}}}
\def\texshopboxnb#1{\boxtextnb{462pt}{\vskip-1.5pc\pshade{\vskip-1.0pc#1\vskip-2.5pc}}}


\input miniltx

\ifx\pdfoutput\undefined
  \def\Gin@driver{dvips.def} 
\else
  \def\Gin@driver{pdftex.def} 
\fi

\input graphicx.sty
\resetatcatcode

\long\def\fig#1#2#3{\vbox{\vskip1pc\vskip#1
\prevdepth=12pt \baselineskip=12pt
\vskip1pc
\hbox to\hsize{\hfill\vtop{\hsize=30pc\noindent{\eightbf Figure #2\ }
{\eightpoint#3}}\hfill}}}

\def\show#1{}

\rightheadline{\botmark}

\pageno=1

\rightheadline{\botmark}

\pn {\bf May 2015 (Revised August 2016)}\hfill{\bf  Report LIDS-P-3173}
\bigskip \bigskip\bigskip

\bigskip

\def\longpapertitle#1#2#3{{\bf \centerline{\helbigb
{#1}}}\medskip{\bf \centerline{\helbigb
{#2}}}\bigskip{\bf \centerline{
{#3}}}\bigskip}

\vskip-3pc

\def\jstar{J^{\raise0.04pt\hbox{\sevenpoint *}} }
\def\qstar{Q^{\raise0.04pt\hbox{\sevenpoint *}} }

\longpapertitle{Regular Policies in}{Abstract Dynamic Programming}{ {Dimitri P.\ Bertsekas\footnote{\dag}{\ninepoint  Dimitri Bertsekas is with the Dept.\ of Electr.\ Engineering and
Comp.\ Science, and the Laboratory for Information and Decision Systems, M.I.T., Cambridge, Mass., 02139.  Many helpful discussions with Huizhen (Janey) Yu on the subject of this paper are gratefully acknowledged.}}}

\vskip-0.5pc
\centerline{\bf Abstract}

We consider challenging dynamic programming models where the associated Bellman equation, and the value and policy iteration algorithms commonly exhibit complex and even pathological behavior. Our analysis is based on the new notion of regular policies. These are policies that are well-behaved with respect to value and policy iteration, and are patterned after proper policies, which are central in the theory of stochastic shortest path problems. We show that the optimal cost function over regular policies may have favorable value and policy iteration properties, which the optimal cost function over all policies need not have. We accordingly develop a unifying methodology to address long standing analytical and algorithmic  issues in broad classes of undiscounted models, including stochastic and minimax shortest path problems, as well as positive cost, negative cost, risk-sensitive, and multiplicative cost problems.

\vskip-1.5pc

\section{Introduction}

\vskip-0.5pc

\pn The purpose of this paper is to address complicating issues that relate to the solutions of Bellman's equation, and the convergence of the value and policy iteration algorithms in total cost infinite horizon dynamic programming (DP for short). We do this in the context of abstract DP, which aims to unify the analysis of DP models and to highlight their fundamental structures. 

To describe broadly our analysis, let us note two types of models. The first is the {\it contractive models\/}, introduced in [Den67], which involve an abstract DP mapping that is a contraction over the space of bounded functions over the state space. These models apply primarily in discounted infinite horizon problems  of various types, with bounded cost per stage. The second is the {\it noncontractive models\/}, developed in [Ber75] and [Ber77] (see also [BeS78], Ch.\ 5), for which the abstract DP mapping is not a contraction of any kind but is instead monotone. Among others, these models include shortest path problems of various types, as well as the classical nonpositive and nonnegative cost DP problems, introduced in [Bla65] and [Str66], respectively. It is well known that contractive models are analytically and computationally well-behaved, while noncontractive models exhibit significant pathologies, which interfere with their effective solution.

In this paper we focus on {\it semicontractive models\/} that were introduced in the recent monograph [Ber13]. These models are characterized by an abstract DP mapping, which for some policies has a contraction-like property, while for others it does not. A central notion in this regard is   {\it $S$-regularity\/} of a stationary policy, where $S$ is a set of cost functions. This property, defined formally in Section 5, is related to classical notions of asymptotic stability, and it roughly means that value iteration using that policy converges to the same limit, the cost function of the policy, for every starting function in the set $S$. 

A prominent case where regularity concepts are central is finite-state  problems of finding an optimal stochastic shortest path (SSP for short). These are Markovian decision problems involving a termination state, where one aims to drive the state of a Markov chain to a termination state at minimum expected cost. They have been discussed in many sources, including the books [Pal67], [Der70], [Whi82], [Ber87], [BeT89],  [BeT91], [Put94], [HeL99], and [Ber12], where they are sometimes referred to by earlier names such as ``first passage problems" and ``transient programming problems." Here some stationary policies called  {\it proper} are guaranteed to terminate starting from every initial state, while others called  {\it improper} are not. The proper policies involve a (weighted sup-norm) contraction mapping and are $S$-regular (with $S$ being the set of real-valued functions over the state space), while the improper ones are not.

The notion of $S$-regularity of a stationary policy is patterned after the notion of a proper policy, but applies more generally in abstract DP. It was used extensively in [Ber13], and in the subsequent papers [Ber15a] and [Ber16] as a unifying analytical vehicle for a variety of total cost stochastic and minimax problems. A key idea is that the optimal cost function over $S$-regular policies only, call it $\jstar_S$, is the one produced by the standard algorithms,  starting from functions $J\in S$ with $J\ge \jstar_S$. These are the value and policy iteration algorithms (abbreviated as VI and PI, respectively), as well as algorithms based on linear programming and related methods. By contrast, the optimal cost function over all policies $\jstar$ may not be obtainable by these algorithms, and indeed $\jstar$ may not be a solution of Bellman's equation; this can happen in particular in SSP problems with zero length cycles (see an example due to [BeY16], which also applies to multiplicative cost problems [Ber16]). 

One purpose of this paper is to extend the notion of $S$-regularity to nonstationary policies, and to demonstrate the use of this extension for establishing convergence of  VI and PI. We show that for important special cases of optimal control problems, our approach yields substantial improvements over the current state of the art, and highlights the fundamental convergence mechanism of  VI and PI  in semicontractive models. A second purpose of the paper is to use the insights of the nonstationary policies extension to refine the stationary regular policies analysis of [Ber13], based on PI-related properties of the set $S$. The paper focuses on issues of  existence and uniqueness of solution of Bellman's equation, and the convergence properties of the VI and PI  algorithms, well beyond the analysis of [Ber13]. A more extensive treatment of the subject of the paper (over 100 pages), which includes elaborations of the analysis, examples, and applications, is given in unpublished internet-posted updated versions of Chapters 3 and 4 of [Ber13], which may be found in the author's web site 
(http://web.mit.edu/dimitrib/www/abstractdp\_MIT.html). 

The paper is organized as follows. After  formulating our abstract DP model in Section 2, we develop the main ideas of the regularity approach for nonstationary policies in Section 3. In Section 4 we illustrate our results by applying them to  nonnegative cost stochastic optimal control problems, and we discuss the convergence of VI, following the analysis of the paper  [YuB13]. In Sections 5-7, we specialize the notion of $S$-regularity to stationary policies, and we refine and streamline the analysis given in the monograph [Ber13], Chapter 3. As an example, we establish the convergence of VI and PI under new and easily verifiable conditions in undiscounted deterministic optimal control problems with a terminal set of states. Other applications of the theory of Sections 5-7 are given in [Ber15a] for robust (i.e., minimax) shortest path planning problems, and in [Ber16] for the class of affine monotonic models, which includes multiplicative and risk sensitive/exponential cost models. 

\vskip-1pc

\section{Abstract Dynamic Programming Model}
\vskip-0.5pc

\pn  We review the abstract DP model that will be used throughout this paper (see Section 3.1 of [Ber13]). Let $X$  and $U$ be two sets, which we  refer to as a set of ``states" and a set of ``controls," respectively. For each $x\in X$, let $U(x)\subset U$ be a nonempty subset of controls that are feasible at state $x$. We denote by ${\cal M}$ the set of all functions $\m:X\mapsto U$ with $\m(x)\in U(x)$, for all $x\in X$. 

We consider policies, which are sequences $\p=\{\m_0,\m_1,\ldots\}$, with $\m_k\in{\cal M}$ for all $k$. We denote by $\Pi$ the set of all policies. We refer to a sequence $\{\m,\m,\ldots\}$, with $\m\in{\cal M}$,  as  a {\it stationary policy\/}. With slight abuse of terminology, we will also refer to any $\m\in{\cal M}$ as a ``policy" and use it in place of $\{\m,\m,\ldots\}$, when confusion cannot arise.

We denote by $\re$ the set of real numbers, by $R(X)$ the set of real-valued functions $J:X\mapsto\re$, and by $E(X)$  the subset of extended real-valued functions $J:X\mapsto\re\cup\{-\infty,\infty\}$. We denote by $E^+(X)$ the set of all nonnegative extended real-valued functions of $x\in X$. Throughout the paper, when we write $\lim$, $\limsup$, or $\liminf$ of a sequence of functions we mean it to be pointwise. We also write $J_k\to J$ to mean that $J_k(x)\to J(x)$ for each $x\in X$, and we write $J_k\downarrow J$ if $\{J_k\}$ is monotonically nonincreasing and $J_k\to J$.

We introduce a mapping $H:X\times U\times E(X)\mapsto\re\cup\{-\infty,\infty\}$, 
satisfying the following condition.

\xdef\assumptionmon{\assumptionn}\assumptionnum\show{myproposition}

\texshopbox{\assumption{\assumptionmon: (Monotonicity)} 
If $J,J'\in E(X)$ and $J\le J'$, then
$$H(x,u,J)\le H(x,u,J'),\qquad \forall\ x\in X,\ u\in U(x).$$
\vskip-1pc}

We define the mapping $T$ that maps a function $J\in E(X)$ to the function $TJ\in E(X)$, given by
$$(TJ)(x)=\inf_{u\in U(x)}H(x,u,J),\qquad \forall\ x\in X,\, J\in E(X).\old{\eqnum\show{oneo}}$$
Also for each $\m\in{\cal M}$, we define the mapping $T_\m:E(X)\mapsto E(X)$ by
$$(T_\m J)(x)=H\big(x,\m(x),J\big),\qquad \forall\ x\in X,\, J\in E(X).\old{\eqnum\show{oneo}}$$
The monotonicity assumption implies the following properties for all $J,J'\in E(X)$, and $k=0,1,\ldots$,
$$J\le J'\qquad\implies\qquad T^kJ\le T^kJ',\qquad T_\m^kJ\le T_\m^kJ',\quad \forall\ \m\in{\cal M},$$
$$J\le TJ\qquad\implies\qquad T^kJ\le T^{k+1}J,\qquad T_\m^kJ\le T_\m^{k+1}J,\quad \forall\ \m\in{\cal M},$$
which will be used repeatedly in what follows. Here $T^k$ and $T_\m^k$ denotes the composition of $T$ and $T_\m$, respectively, with itself $k$ times. More generally, given $\m_0,\ldots,\m_k\in{\cal M}$, we denote by $T_{\m_0}\cdots T_{\m_k}$ the composition of $T_{\m_0},\ldots,T_{\m_k}$, so for all $J\in E(X)$, 
$$(T_{\m_0}\cdots T_{\m_k}J\big)(x)=\big(T_{\m_0}\big(T_{\m_1}\cdots \big(T_{\m_{k-1}}(T_{\m_k}J)\big)\cdots\big)\big)(x),\qquad\forall\ x\in X.$$

We next consider cost functions associated with $T_\m$ and $T$.  We introduce a function $\bar J\in E(X)$, and we define the infinite horizon cost of a policy as the upper limit of its finite horizon costs with $\bar J$ being the cost function at the end of the horizon (limit cannot be used since it may not exist).

\xdef\definitioncost{\defn}\defnum\show{myproposition}

\texshopbox{\definition{\definitioncost:}Given a function $\bar J\in E(X)$, for a policy $\p\in \Pi$ with $\p=\{\m_0,\m_1,\ldots\}$, we define the cost function of $\p$ by
$$J_\p(x)=\limsup_{k\to\infty}\,(T_{\m_0}\cdots T_{\m_k} \bar J)(x),\qquad \forall\ x\in X.\xdef\polcost{\lab}\eqnum\show{oneo}$$
The optimal cost function $\jstar$ is defined by
$$\jstar(x)=\inf_{\p\in\P}J_\p(x),\qquad\forall\ x\in X.$$
A policy $\p^*\in\P$ is said to be optimal if $J_{\p^*}=\jstar$.  
}

The model just described is broadly applicable, and includes as special cases nearly all the interesting types of total cost infinite horizon DP problems, including stochastic and minimax, discounted and undiscounted, semi-Markov, multiplicative, risk-sensitive, etc (see [Ber13]).\footnote{\dag}{\ninepoint However, our model cannot address those stochastic DP models where measurability issues are an important mathematical concern. In the stochastic optimal control problem of Example 2.1, we bypass these issues by assuming that the disturbance space is countable, which includes the deterministic system case, and the case where the system is stochastic with a countable state space (e.g., a countable state Markovian decision problem). Then, the expected value needed to express the finite horizon cost of a policy [cf.\ Eq.\ \polcost] can be written as a summation over a countable index set, and is well-defined  for all policies, measurable or not.} The following is a stochastic optimal control problem, which we will use in this paper both to obtain new results and also as a vehicle to  illustrate our approach.

\xdef\exampledisc{\exampl}\examplnum\show{myexample}

\beginexample{\exampledisc\ (Stochastic Optimal Control - Undiscounted Markovian Decision Problems)}Consider an infinite horizon stochastic optimal control problem involving a stationary discrete-time dynamic system
where the state
is an element of a space $X$, and the control is an element of a space $U$.  
 The control $u_k$ is constrained to take values in a given nonempty subset $U(x_k)$
of $U$, which depends on the current state $x_k$ [$u_k
\in U(x_k)$, for all $x_k\in X$]. For a policy $\p=\{\m_0,\m_1,\ldots\}$, the state evolves according to a system equation
$$x_{k+1} = f\big(x_k,\m_k(x_k),w_k\big),\qquad k =
0,1,\ldots,\xdef\systemeq{\lab}\eqnum\show{oneo}$$ 
where $w_k$ is a random disturbance that takes values from a space $W$. We assume that $w_k$, $k=0,1,\ldots$, are characterized
by probability distributions $P(\cdot\mid x_k,u_k)$ that are identical for all $k$, where
$P(w_k\mid x_k,u_k)$ is the probability of occurrence of $w_k$, when the current state and control
are $x_k$ and $u_k$, respectively.  Thus the probability of $w_k$ may depend explicitly on $x_k$ and
$u_k$, but not on values of prior disturbances $w_{k-1},\ldots,w_0$. We allow infinite state and control spaces,  as well as problems with discrete (finite or countable) state space (in which case the underlying system is a Markov chain). However, for technical reasons that relate to measure theoretic issues, we assume that $W$ is a countable set. 

Given an initial state $x_0$, we want to find a policy $\pi = \{\mu_0,\mu_1,\ldots\}$,
where $\mu_k:X \mapsto U$, $\mu_k(x_k) \in U(x_k)$, for all $x_k \in X$, $k=0,1,\ldots$, that
minimizes
$$J_\pi(x_0) = \limsup_{k\tends\infty} E\lf\{\sum_{ t=0}^{k}
\a^k g\bl(x_ t,\mu_ t(x_ t), w_ t\br)\ri\},$$
subject to the system equation constraint \systemeq,
where $g$ is the one-stage cost function, and $\a\in(0,1]$ is the discount factor. This is a classical problem, which is discussed extensively in various sources, such as the books [BeS78], [Whi82], [Put94], [Ber12]. Under very mild conditions guaranteeing that Fubini's theorem can be applied (see [BeS78], Section 2.3.2), it coincides with the abstract DP problem that corresponds to the mapping 
$$H(x,u,J)=E\big\{g(x,u,w)+\a J\big(f(x,u,w)\big)\big\},\xdef\socmapdisc{\lab}\eqnum\show{oneo}$$ 
and $\skew5\bar J(x)\equiv0$. Here, $(T_{\m_0}\cdots T_{\m_{k}}\skew6\bar J)(x)$ is the expected cost of the first $k+1$ periods using $\p$ starting from $x$, and with terminal cost 0 (the value of $\skew6\bar J$ at the terminal state).

\endexample

\vskip-2pc

\section{Regular Policies, Value Iteration, and Fixed Points of $T$}

\pn Generally, in an abstract DP model, one expects to establish that $\jstar$ is a fixed point of $T$. This is  known to be true for most DP models under reasonable conditions, and in fact it may be viewed as an indication of exceptional behavior when it does not hold. The fixed point equation $J=TJ$, in the context of standard special cases, is the classical {\it Bellman equation\/}, the centerpiece of infinite horizon DP.  For some abstract DP models, $\jstar$ is the unique fixed point of $T$ within a convenient subset of $E(X)$; for example, contractive models where $T_\m$ is a contraction mapping for all $\m\in{\cal M}$, with respect to some norm and with a common modulus of contraction. However, in general $T$ may have multiple fixed points within $E(X)$, including for some popular DP problems, while in exceptional cases, $\jstar$ may not be among the fixed points of $T$ (see [BeY16] for a relatively simple SSP example of this type). 

A related question is the convergence of VI. This is the algorithm that generates $T^kJ$, $k=0,1,\ldots,$ starting from a function $J\in E(X)$. Generally, for abstract DP models where $\jstar$ is a fixed point of $T$, VI converges to $\jstar$ starting from within some subset of initial functions $J$, but not from every $J$; this is certainly true when $T$ has multiple fixed points.
One of the  purposes of this paper is to characterize the set of functions starting from which VI converges to $\jstar$,  and the related issue of multiplicity of fixed points, through notions of regularity that we  now introduce.

\xdef\definitioncregular{\defn}\defnum\show{myproposition}

\texshopbox{\definition{\definitioncregular:}For a nonempty set of functions $S\subset E(X)$, we say that a set ${\cal C}$ of policy-state pairs $(\p,x)$, with $\p\in \Pi$ and $x\in X$,  is  {\it $S$-regular} if
$$J_\p(x)=\limsup_{k\to\infty}\,(T_{\m_0}\cdots T_{\m_k}J)(x),\qquad \forall\ (\p,x)\in {\cal C},\ J\in S.
$$
}

A nonempty set ${\cal C}$ of policy-state pairs $(\p,x)$ may be $S$-regular for many different sets $S$. The largest such set is
$$S_{\cal C}=\lf\{J\in E(X)\ \Big|\   J_\p(x)=\limsup_{k\to\infty}\,(T_{\m_0}\cdots T_{\m_k}J)(x),\, \forall\ (\p,x)\in {\cal C}\ri\},$$
and for any nonempty $S\subset S_{\cal C}$, we have that ${\cal C}$ is $S$-regular. Moreover, the set $S_{\cal C}$ is nonempty, since it contains $\bar J$.
For a given ${\cal C}$, consider the function $\jstar_{\cal C}\in E(X)$, given by
$$\jstar_{\cal C}(x)=\inf_{ \{\p\,\mid\, (\p,x)\in{\cal C}\}}J_\p(x),\qquad x\in X.$$
Note that $\jstar_{\cal C}(x)\ge \jstar(x)$ for all $x\in X$ [for those $x\in X$ for which the set of policies $\{\p\,\mid\, (\p,x)\in{\cal C}\}$ is empty, we have $\jstar_{\cal C}(x)=\infty$]. 
We will try to characterize the sets of fixed points of $T$ and limit points of VI in terms of the function  $\jstar_{\cal C}$ for an $S$-regular set ${\cal C}$. The following is a key proposition. In this proposition as well as later {\it when referring to a set ${\cal C}$ that is $S$-regular, we implicitly assume that ${\cal C}$ and $S$ are nonempty\/}.

\xdef\propregset{\propn}\propnum\show{myproposition}

\texshopbox{\proposition{\propregset:} Given a set $S\subset E(X)$, let ${\cal C}$ be an $S$-regular set.
\nitem{(a)} For all $J\in S$, we have
$$\liminf_{k \to\infty}T^k J\le \limsup_{k \to\infty}T^k J\le \jstar_{\cal C}.$$
\nitem{(b)} For all $J'\in E(X)$ with $J'\le TJ'$, and all $J\in E(X)$ such that $J'\le J\le \tl J$ for some $\tl J\in S$, we have
$$J'\le \liminf_{k \to\infty}T^k J\le \limsup_{k \to\infty}T^k J\le \jstar_{\cal C}.$$
}

\proof (a) Using the generic relation $TJ\le T_\m J$, $\m\in{\cal M}$,  and the monotonicity of $T$ and $T_\m$, we have for all $k$
$$(T^k J)(x)\le (T_{\m_0}\cdots T_{\m_{k-1}} J)(x),\qquad \forall\ (\p,x)\in  {\cal C},\ J\in S.$$
By letting $k\to\infty$ and by using the definition of $S$-regularity, it 
follows that
$$\liminf_{k \to\infty}(T^k J)(x)\le \limsup_{k \to\infty}(T^k J)(x)\le \limsup_{k \to\infty}(T_{\m_0}\cdots T_{\m_{k-1}} J)(x)= J_\p(x),\qquad \forall\ {(\p,x)}\in {\cal C},\ J\in S,$$
and taking infimum of the right side over $\big\{\p\mid (\p,x)\in{\cal C}\big\}$, we obtain the result.
\smskip
\pn (b) Using the hypotheses $J'\le TJ'$, and $J'\le J\le \tl J$ for some $\tl J\in S$, and the monotonicity of $T$, we have 
$$J'(x)\le (TJ')(x)\le\cdots \le(T^kJ')(x)\le (T^k J)(x)\le (T^k \tl J)(x).$$
Letting $k\to\infty$ and using part (a), 
we obtain the result. 
\qed
 
\xdef \figregular{\figr}\figrnum\show{myfigure}

Part (b) of the proposition shows that given a set $S\subset E(X)$, a nonempty set ${\cal C}\subset \P\times X$ that is $S$-regular, and a function $J'\in E(X)$ with $J'\le TJ'\le \jstar_{\cal C}$, the convergence of VI is characterized by the {\it valid start region} 
$$\big\{J\in E(X)\mid J'\le J\le \tl J \hbox{ for some }\tl J\in S\big\},$$
and the {\it limit region} 
$$\big\{J\in E(X)\mid J'\le J\le \jstar_{\cal C}\big\}.$$
 The VI algorithm, starting from the former, ends up asymptotically within the latter; cf.\ Fig.\ \figregular.  Note that 
both of these regions depend on ${\cal C}$ and $J'$.

\midinsert
\centerline{\hskip0pc\includegraphics[width=4.1in]{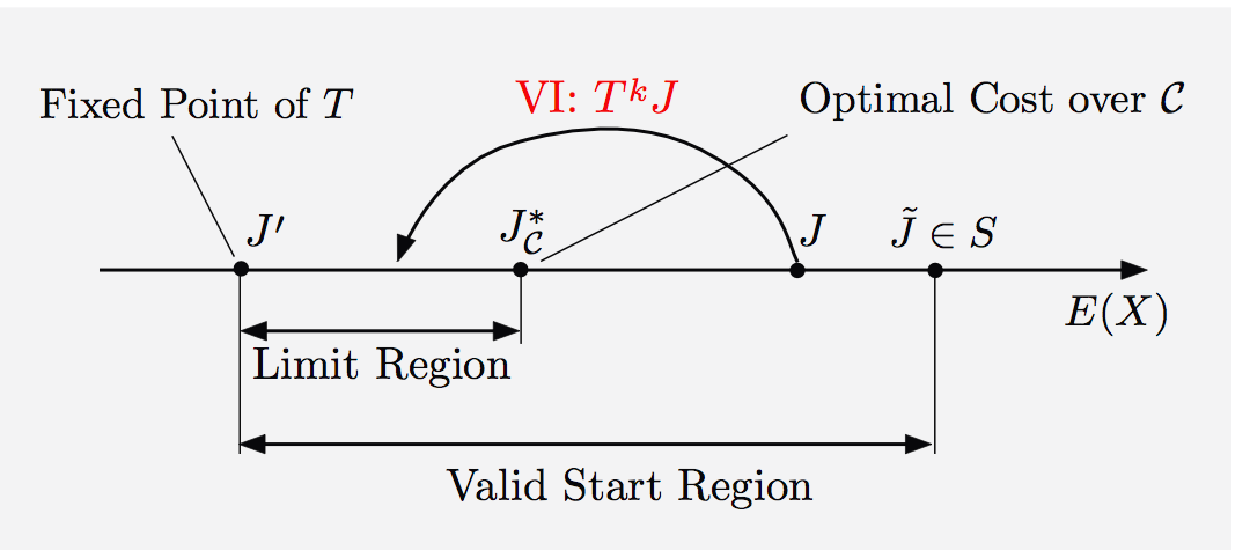}}
\vskip-1pc
\hskip-4pc\fig{0pc}{\figregular.} {Illustration of Prop.\ \propregset. Neither $J^*_{\cal C}$ nor $J^*$ need to be fixed points of $T$, but if ${\cal C}$ is $S$-regular, and there exists $\skew5\tl J\in S$ with $J^*_{\cal C}\le \skew5\tl J$, then $J^*_{\cal C}$ demarcates from above the range of fixed points of $T$ that lie below $\skew5\tl J$.}
\endinsert

 The significance of the preceding property depends of course on the choice of ${\cal C}$ and $S$. With an appropriate choice, however, there are important implications regarding the location of the fixed points of $T$ and the convergence of VI from a broad range of starting points. Some of these implications are the following:

\nitem{(a)} $\jstar_{\cal C}$ is an upper bound to every fixed point $J'$ of $T$ that lies below some $\tl J\in S$ (i.e., $J'\le \tl J$).

\nitem{(b)} If $\jstar_{\cal C}$ is a fixed point of $T$ (an important case for our subsequent development), then VI converges to $\jstar_{\cal C}$ starting from any $J\in E(X)$ such that $\jstar_{\cal C}\le J\le \tl J$ for some $\tl J\in S$. For future reference, we state this result as a proposition.

\xdef\propregsetcorth{\propn}\propnum\show{myproposition}

\texshopbox{\proposition{\propregsetcorth:}  Given a set $S\subset E(X)$, let ${\cal C}$ be an $S$-regular set  and assume that
$\jstar_{\cal C}$ is a fixed point of $T$. 
 Then $\jstar_{\cal C}$ is the only possible fixed point of $T$ within the set of all $J\in E(X)$ such that $\jstar_{\cal C}\le J\le \tl J$ for some $\tl J\in S$. Moreover, $T^kJ\to\jstar_{\cal C}$ for all $J\in E(X)$ such that $\jstar_{\cal C}\le J\le \tl J$ for some $\tl J\in S$.
}

\proof Let $J\in E(x)$ and $\tl J\in S$ be such that $\jstar_{\cal C}\le J\le \tl J$. Using the fixed point property of $\jstar_{\cal C}$ and the monotonicity of $T$, we have
$$\jstar_{\cal C}=T^k\jstar_{\cal C}\le T^kJ\le T^k\tl J,\qquad k=0,1,\ldots.$$
From Prop.\ \propregset(b), with $J'=\jstar_{\cal C}$, it follows that $T^k\tl J\to \jstar_{\cal C}$, so taking limit in the above relation as $k\to\infty$, we obtain $T^kJ\to \jstar_{\cal C}$.
\qed 

The preceding proposition takes special significance when ${\cal C}$ is rich enough so that $\jstar_{\cal C}=\jstar$, as for example in the case where
${\cal C}$ is the set $\Pi\times X$ of all $(\p,x)$,
or other choices to be discussed later.
It then follows that VI converges to $\jstar$ starting from any $J\in E(X)$ such that $\jstar\le J\le \tl J$ for some $\tl J\in S$.\footnote{\dag}{\ninepoint For this statement to be meaningful, the set $\big\{\skew5\tl J\in E(X)\mid J^*\le \skew5\tl J\big\}$ must be nonempty. Generally, it is possible that this set is empty, even though $S$ is assumed nonempty.} In the particular applications to be discussed in Section 4 we will use such a choice. 

Note that Prop.\ \propregsetcorth\ does not say anything about fixed points of $T$ that lie below $\jstar_{\cal C}$. In particular, it does not address the question whether $\jstar$ is a fixed point of $T$, or whether VI converges to $\jstar$ starting from $\bar J$ or from below $\jstar$; these are major questions in abstract DP models, which are typically handled by special analytical techniques that are tailored to the particular model's structure and assumptions. Significantly, however,  these questions have been already answered in the context of various  models, and when available, they can be used to supplement the preceding propositions. For example, the DP books  [Pal67], [Der70], [Whi82], [Put94], [HeL99], [Ber12], [Ber13] provide extensive analysis for the most common infinite horizon stochastic optimal control problems: discounted, SSP, nonpositive cost, and nonnegative cost problems. 

In particular, for discounted problems [the case of the mapping \socmapdisc\ with $\a\in(0,1)$ and $g$ being a bounded function], underlying sup-norm contraction properties guarantee that $\jstar$ is the unique fixed point of $T$ within the class of bounded real-valued functions over $X$, and that VI converges to $\jstar$ starting from within that class. This is also true for finite-state SSP problems, involving a cost-free termination state, under some favorable conditions (there must exist a proper policy, i.e., a stationary policy that leads to the termination state with probability 1, improper policies must have infinite cost for some states, and some finiteness or compactness conditions on the control space $U$ must be satisfied; see [BeT91], [Ber12]). 

The  paper [BeY16] also considers finite-state SSP problems, but under the weaker assumptions that there exists at least one proper policy, that $\jstar$ is real-valued, and $U$ satisfies some finiteness or compactness conditions. Under these assumptions, $\jstar$ need not be a fixed point of $T$, as shown in [BeY16] with an example. In the context of the present paper, a useful choice is to take 
${\cal C}=\big\{(\m,x)\mid \m\,\hbox{: proper}\big\},$
 in which case $\jstar_{\cal C}$ is the optimal cost function that can be achieved using proper policies only. It was shown in [BeY16] that $\jstar_{\cal C}$ is a fixed point of $T$, so by Prop.\ \propregsetcorth, VI converges to $\jstar_{\cal C}$ starting from any real-valued $J\ge \jstar_{\cal C}$. 
  
For nonpositive and nonnegative cost problems (cf.\ Example \exampledisc\ with $g\le0$ or $g\ge0$, respectively), $\jstar$ is a fixed point of $T$, but not necessarily unique.
However, for nonnegative cost problems, some new results on the existence of fixed points of $T$ and convergence of VI were recently proved in [YuB13]. It turns out that one may prove these results by using Prop.\ \propregsetcorth, with 
an appropriate choice of ${\cal C}$.
The proof  uses the arguments of Appendix E of  [YuB13], and will  be given in Section 4.1.

A class of DP problems with more complicated structure is the general convergence model discussed in the thesis [Van81] and the survey paper  [Fei02]. This is the case of Example \exampledisc\ where the cost per stage $g$ can take both positive and negative values, under some restrictions that guarantee that $J_\p$ is defined by Eq.\ \polcost\ as a limit. The paper [Yu15] describes the complex issues of convergence of VI for these models, and in an infinite space setting that addresses measurability issues. We note that there are examples of general convergence models where $X$ and $U$ are finite sets, but VI does not converge to $\jstar$ starting from $\bar J$ (see Example 3.2 of [Van81], Example 6.10 of [Fei2], and Example 4.1 of [Yu15]). The analysis of [Yu15] may also be used to bring to bear Prop.\  \propregset\ on the problem, but this analysis is beyond our scope in this paper.

\vskip-0.5pc

\subsubsection{\bf The Case Where $\jstar_{\cal C}\le \bar J$}

\pn It is well known that the results for nonnegative cost and nonpositive cost infinite horizon stochastic optimal control problems are markedly different. In particular, roughly speaking, PI behaves better when the cost is nonnegative, while VI behaves better if the cost is nonpositive. These differences extend to the so-called {\it monotone increasing} and {\it monotone decreasing} abstract DP models, where a principal assumption is that $T_\m\bar J\ge \bar J$ and $T_\m\bar J\le \bar J$ for all $\m\in{\cal M}$, respectively (see [Ber13], Ch. 4). In the context of regularity, with ${\cal C}$ being $S$-regular, it turns out that there are analogous significant differences between the cases $\jstar_C\ge \bar J$ and $\jstar_C\le \bar J$. The following proposition establishes some favorable aspects of the condition $\jstar_{\cal C}\le \bar J$ in the context of VI. These can be attributed to the fact that $\bar J$ can always be added to $S$ without affecting the $S$-regularity of ${\cal C}$, so $\bar J$ can serve as the element $\tl J$ of $S$ with $\jstar_C\le \tl J$ in Props.\ \propregset\ and \propregsetcorth\ (see the proof of the following proposition).

\xdef\propregsetcoro{\propn}\propnum\show{myproposition}
\vskip-0.5pc

\texshopbox{\proposition{\propregsetcoro:}  Given a set $S\subset E(X)$, let ${\cal C}$ be an $S$-regular set and assume that $\jstar_{\cal C}\le \bar J$. Then:
\nitem{(a)} For all $J'\in E(X)$ with $J'\le TJ'$, we have%
$$J'\le \liminf_{k\to\infty}T^k\bar J\le \limsup_{k\to\infty}T^k\bar J\le \jstar_{\cal C}.$$
\nitem{(b)} If $\jstar_{\cal C}$ is a fixed point of $T$, then $\jstar=\jstar_{\cal C}$ and we have $T^k \bar J\to \jstar$ as well as $T^k J\to \jstar$ for every $J\in E(X)$ such that $\jstar\le J\le \tl J$ for some $\tl J\in S$.¡
}

\proof (a) If $S$ does not contain $\bar J$, we can replace $S$ with $\bar S=S\cup\{\bar J\}$, and ${\cal C}$ will still be $\bar S$-regular. By applying Prop.\ \propregset(b) with $S$ replaced by $\bar S$ and $\tl J=\bar J$, the result follows.

\smskip
\pn (b) Assume without loss of generality that $\bar J\in S$ [cf.\ the proof of part (a)]. By using Prop.\ \propregsetcorth\ with $\tl J=\bar J$, we have $\jstar_{\cal C}=\lim_{k\to\infty}T^k\bar J$.
This relation yields for any policy $\p=\{\m_0,\m_1,\ldots\}\in \P$, 
$$\jstar_{\cal C}=\lim_{k\to\infty}T^k\bar J\le \limsup_{k\to\infty}T_{\m_0}\cdots T_{\m_{k-1}}\bar J=J_\p,$$
so by taking the infimum over $\p\in \P$, we obtain 
$\jstar_{\cal C}\le \jstar$. Since generically we have $\jstar_{\cal C}\ge \jstar$, it follows that $\jstar_{\cal C}=\jstar$. Finally, from Prop.\ \propregsetcorth, we obtain $T^k J\to \jstar$ for all $J\in E(X)$ such that $\jstar\le J\le \tl J$ for some $\tl J\in S$. \qed

\vskip-0.5pc

As a special case of the preceding proposition, we have that if $\jstar\le \bar J$ and $\jstar$ is a fixed point of $T$, then $\jstar=\lim_{k\to\infty}T^k\bar J$, and for every other fixed point $J'$ of $T$ we have $J'\le \jstar$ (apply the proposition with ${\cal C}=\Pi\times X$ and $S=\{\bar J\}$, in which case $\jstar_{\cal C}=\jstar\le \bar J¡$). This special case is relevant, among others, to the monotone decreasing models (see [Ber13], Section 4.3), where $T_\m\bar J\le \bar J$ for all $\m\in{\cal M}$, in which case it is known that $\jstar$ is a fixed point of $T$ under mild conditions.  We then obtain a classical result on the convergence of VI  for nonpositive cost models.
The proposition also applies to a classical type of search problem with both positive and negative costs per stage. This is Example \exampledisc, where at each $x\in X$  we have  $E\big\{g(x,u,w)\big\}\ge0$ for all $u$ except one that leads to a termination state with probability 1 and nonpositive cost. 
Note that without the assumption $\jstar_{\cal C}\le \bar J$ in the preceding proposition, it is possible that $T^k\bar J$ does not converge to $\jstar$, even if $\jstar_{\cal C}=\jstar=T\jstar$, as is well known in the theory of nonnegative cost infinite horizon stochastic optimal control. 

Generally, it is important to choose properly the set ${\cal C}$ in order to obtain meaningful results. Note, however, that in a given problem the interesting choices of ${\cal C}$ are usually limited, and that the propositions of this section can guide a favorable choice. One useful approach is to try the set
$${\cal C}=\big\{(\p,x)\mid J_\p(x)<\infty\big\},$$
so that $\jstar_{\cal C}=\jstar$. By the definition of regularity, if $S$ is any subset of the set
$$S_{\cal C}=\lf\{J\in E(X)\ \Big|\   J_\p(x)=\limsup_{k\to\infty}(T_{\m_0}\cdots T_{\m_k}J)(x),\, \forall\ (\p,x)\in {\cal C}\ri\},$$
then ${\cal C}$ is $S$-regular. One may then try to derive a suitable subset of $S_{\cal C}$ that admits an interesting characterization. This is the approach followed in the applications of the next section. Another approach, discussed in Section 5, is to focus on an interesting subset $\ol {\cal M}$ of {\it stationary} policies such that for the set
$${\cal C}=\ol {\cal M}\times X,$$
we have $\jstar_{\cal C}=\jstar$.

\vskip-1.5pc

\section{Applications in Stochastic Optimal Control}

\vskip-0.5pc

\pn In this section, we will consider the stochastic optimal control problem of  Example \exampledisc, where
$$H(x,u,J)=E\big\{g(x,u,w)+\a J\big(f(x,u,w)\big)\big\},\xdef\socmapundiscdisc{\lab}\eqnum\show{oneo}$$ 
and $\bar J(x)\equiv0$. Here $\a\in(0,1]$ is the discount factor and we assume that the expected cost per stage is nonnegative:
$$0\le E\big\{g(x,u,w)\big\}< \infty,\qquad\forall\ x\in X,\ u\in U(x).\xdef\finiteexpg{\lab}\eqnum\show{oneo}$$
This is a classical problem, also known as the negative DP model [Str66].

We will use some classical results for this problem, which we collect in the following proposition (for proofs, see e.g., [BeS78], Props.\ 5.2, 5.4, and 5.10, or [Ber13], Props.\ 4.3.3, 4.3.9, and 4.3.14).

\vskip-0.5pc
\xdef\propnegdp{\propn}\propnum\show{myproposition}

\texshopbox{\proposition{\propnegdp:} Consider the stochastic optimal control problem where $H$ is given by Eq.\ \socmapundiscdisc, $g$ satisfies the nonnegativity condition \finiteexpg, and $\a\in(0,1]$. Then: 
\nitem{(a)} $\jstar=T\jstar$ and  if $J\in E^+(X)$ satisfies $J\ge TJ$, then  $J\ge\jstar$.
\nitem{(b)} For all $\m\in{\cal M}$ we have $J_\m=T_\m J_\m$.
\nitem{(c)} $\m^*\in{\cal M}$ is optimal if and only if $T_{\m^*}\jstar=T\jstar$.
\nitem{(d)} If $U$ is a metric space and the sets 
$$U_k(x,\l) =\big\{ u\in U(x)\mid H(x,u,T^k\bar J)\le\l\big\}\xdef\onetwen{\lab}\eqnum\show{oneo}$$
are compact for all $x\in X$, $\l\in \re$, and $k$,  then there exists at least one optimal stationary policy, and we have $T^kJ\to\jstar$ for all $J\in E^+(X)$ with $J\le \jstar$.}

Note that there may exist fixed points $J'$ of $T$ with $J'\ge \jstar$, while VI or PI may not converge to $\jstar$ starting from above $\jstar$. However, convergence of VI to $\jstar$ from above, if it occurs, is often much faster than convergence from below, so starting points $J\ge \jstar$ may be desirable. One well-known such case is deterministic finite-state  shortest path problems where major algorithms, such as the Bellman-Ford method or other label correcting methods have polynomial complexity, when started from $J$ above $\jstar$, but only pseudopolynomial complexity when started from other initial conditions.

We will now establish conditions for the uniqueness of $\jstar$ as a fixed point of $T$, and the convergence of VI and PI. We will consider separately the cases $\a=1$ and $\a<1$. Our analysis will proceed as follows:
\nitem{(a)} Define a set ${\cal C}$ such that $\jstar_{\cal C}=\jstar$.
\nitem{(b)} Define a set $S\subset E^+(X)$ such that $\jstar\in S$ and ${\cal C}$ is $S$-regular.
\nitem{(c)} Use Prop.\ \propregsetcorth\ in conjunction with the fixed point properties of $\jstar$ [cf.\ Prop.\ \propnegdp(a)] to show that $\jstar$ is the unique fixed point of $T$ within $S$, and that the VI algorithm converges to $\jstar$ starting from $J$ within the set $\{J\in S\mid J\ge \jstar\}$.
\nitem{(d)} Use the compactness condition of Prop.\ \propnegdp(d), to enlarge the set of functions starting from which VI converges to $\jstar$.

\subsection{Nonnegative Undiscounted Cost Stochastic DP}

\pn Assume that the problem is undiscounted, i.e., $\a=1$. Consider the set
$${\cal C}=\big\{(\p,x)\mid J_\p(x)<\infty\big\},$$
for which we have $\jstar_{\cal C}=\jstar$, and assume that ${\cal C}$ is nonempty. 

Let us denote by  $E^\p_{x_0}\{\cdot\}$ the expected value with respect to the probability measure induced by $\p\in \Pi$ under initial state $x_0$, and let us consider the set
$$S=\big\{J\in E^+(X)\mid E_{x_0}^{\p}\big\{J(x_k)\big\}\to0,\ \forall\ {(\p,x_0)}\in {\cal C}\big\}.\xdef\regsets{\lab}\eqnum\show{oneo}$$
We will show that $\jstar\in S$ and that ${\cal C}$ is $S$-regular.
Once this is done, it will follow from Prop.\ \propregsetcorth\ and the fixed point property of $\jstar$ [cf.\ Prop.\ \propnegdp(a)] that $T^k J\to \jstar$ for all $J\in S$ that satisfy $J\ge \jstar$. If the sets 
$U_k(x,\l)$ of Eq.\ \onetwen\ are compact, the convergence of VI starting from below $\jstar$ will also be guaranteed. We have the following proposition. The proof uses the line of argument of Appendix E of  [YuB13].

\xdef\propdesoc{\propn}\propnum\show{myproposition}

\texshopbox{\proposition{\propdesoc:  (Convergence of VI)} Consider the stochastic optimal control problem of this section, assuming $\a=1$ and the cost nonnegativity condition \finiteexpg.
Then $\jstar$ is the unique fixed point of $T$ within $S$, and we have  $T^kJ\to\jstar$ for all $J\ge \jstar$ with $J\in S$. If in addition $U$ is a metric space, and the sets $U_k(x,\l)$ of 
Eq.\ \onetwen\ 
are compact for all $x\in X$, $\l\in \re$, and $k$, we have  $T^kJ\to \jstar$ for all $J\in S$, and an optimal stationary policy is guaranteed to exist.
}

\proof  We have for all $J\in E(X)$, ${(\p,x_0)}\in{\cal C}$, and $k$,
$$(T_{\m_0}\cdots T_{\m_{k-1}}J)(x_0)=E_{x_0}^{\p}\big\{J(x_k)\big\}+E_{x_0}^{\p}\lf\{\sum_{ t=0}^{k-1}g\big(x_ t,\m_ t(x_ t),w_ t\big)\ri\},\xdef\finitehorcost{\lab}\eqnum\show{oneo}$$
where $\m_ t$, $ t=0,1,\ldots$, denote generically the components of $\p$.
By the cost nonnegativity condition \finiteexpg, the rightmost term above converges to $J_\p(x_0)$ as $k\to\infty$, so by taking upper limit, we obtain
$$\limsup_{k\to\infty}(T_{\m_0}\cdots T_{\m_{k-1}}J)(x_0)=\limsup_{k\to\infty}E_{x_0}^{\p}\big\{J(x_k)\big\}+J_\p(x_0).$$
Thus in view of the definition \regsets\ of $S$, we see that for all $(\p,x_0)\in{\cal C}$ and $J\in S$, we have
$$\limsup_{k\to\infty}(T_{\m_0}\cdots T_{\m_{k-1}}J)(x_0)=J_\p(x_0),$$
so ${\cal C}$ is $S$-regular.

 We next show that $\jstar\in S$. We have for all $(\p,x_0)\in{\cal C}$ 
$$J_\p(x_0)=E_{x_0}^{\p}\big\{g\big(x_0,\m_0(x_0),w_0\big)\big\}+E_{x_0}^{\p}\big\{J_\p(x_{1})\big\},$$
 and more generally,
$$E_{x_0}^{\p}\big\{J_\p(x_ t)\big\}=E_{x_0}^{\p}\big\{g\big(x_ t,\m_ t(x_ t),w_ t\big)\big\}+E_{x_0}^{\p}\big\{J_\p(x_{ t+1})\big\},\qquad\forall\  t=0,1,\ldots,\xdef\finitejpi{\lab}\eqnum\show{oneo}$$
where $\{x_t\}$ is the sequence generated starting from $x_0$ and using $\p$. Using the defining property $J_\p(x_0)<\infty$ of ${\cal C}$, it follows  that all the terms in the above relations are finite, and in particular
$$E_{x_0}^{\p}\big\{J_\p(x_ t)\big\}<\infty,\qquad \forall\ (\p,x_0)\in {\cal C},\  t=0,1,\ldots.$$
By adding Eq.\ \finitejpi\ for $t=0,\ldots,k-1$, and canceling the finite terms 
$E_{x_0}^{\p}\big\{J_\p(x_ t)\big\}$ for $ t=1,\ldots,k-1$,
$$J_\p(x_0)=E_{x_0}^{\p}\big\{J_\p(x_k)\big\}+\sum_{ t=0}^{k-1}E_{x_0}^{\p}\big\{g\big(x_ t,\m_ t(x_ t),w_ t\big)\big\},\qquad \forall\  (\p,x_0)\in {\cal C},\ k=1,2,\ldots.$$
The rightmost term above tends to $J_\p(x_0)$ as $k\to\infty$, so we obtain 
$E_{x_0}^{\p}\big\{J_\p(x_k)\big\}\to0$ for all $(\p,x_0)\in {\cal C}.$
Since $0\le \jstar\le J_\p$ for all $\p$, it follows that
$$E_{x_0}^{\p}\big\{\jstar(x_k)\big\}\to0,\qquad \forall\ x_0\hbox{ with }\jstar(x_0)<\infty.$$
Thus $\jstar\in S$.

From Prop.\ \propregsetcorth\ it  follows that $\jstar$ is the unique fixed point of $T$ within $\big\{J\in S\mid J\ge \jstar\big\}$.
On the other hand, every fixed point $J\in E^+(X)$ of $T$ satisfies $J\ge \jstar$ by  Prop.\ \propnegdp(a), so $\jstar$ is the unique fixed point of $T$ within $S$. Also  from Prop.\ \propregsetcorth\ we have that the VI sequence $\{T^kJ\}$ converges to $\jstar$ starting from any $J\in S$ with $J\ge \jstar$. 
Finally, for any $J\in S$, let us select $\tl J\in S$ with $\tl J\ge \jstar$ and $\tl J\ge J$, and note that by the monotonicity of $T$, we have
$T^k\bar J\le T^k J\le T^k \tl J.$ 
If we also assume compactness of the sets $U_k(x,\l)$ of Eq.\ \onetwen, then by Prop.\ \propnegdp(d), we have $T^k\bar J\to \jstar$, which together with the convergence $T^k\tl J\to \jstar$ just proved, implies that $T^k J\to \jstar$. 
\qed

A consequence of the preceding proposition is an interesting condition for VI convergence from above, which was first proved in [YuB13]. In particular, since $\jstar\in S$, any $J$ satisfying $\jstar\le J\le c\jstar$ for some $c>0$ belongs to $S$, so we have the following.

\xdef\propbridgecond{\propn}\propnum\show{myproposition}

\texshopbox{\proposition{\propbridgecond: [YuB13]} We have $T^kJ\to\jstar$ for all $J\in E(X)$ satisfying
$\jstar\le J\le c\jstar$ for some $c>0.$
}

The preceding proposition highlights a requirement for the reliable implementation of VI: it is important to know the sets $X_s=\big\{x\in X\mid \jstar(x)=0\big\}$ and  $X_\infty=\big\{x\in X\mid \jstar(x)=\infty\big\}$ in order to obtain a suitable initial condition $J\in E(X)$ satisfying
$\jstar\le J\le c\jstar$ for some $c>0$. For finite state and control problems, the set $X_s$ can be computed in polynomial time as shown in the paper [BeY16], which also provides a method for dealing with cases where $X_\infty$ is nonempty, based on adding a high cost artificial control st each state. 

Regarding PI, we note that the analysis of Section 5.2 will guarantee its convergence for the stochastic problem of this section if somehow it can be shown that $\jstar$ is the unique fixed point of $T$ within a subset of $\{J\mid J\ge \jstar\}$ that contains the limit $J_\infty$ of PI. This result was given as Corollary 5.2 in [YuB13]. Alternatively, there is a mixed VI and PI algorithm proposed in [YuB13], which can be applied under the condition of Prop.\ \propbridgecond, and applies to a more general problem where $w$ can take an uncountable number of values and measurability issues are an important concern. 

Finally, we note that in this section we do not consider any special structure, other than the expected cost nonnegativity condition \finiteexpg. In particular, we do not discuss the implications of the possible existence of a termination state as in finite-state or countable-state SSP problems. The approach of this paper is relevant to the convergence analysis of VI and PI for such problems, and for a corresponding analysis for finite-state problems, we refer to the paper [BeY16].

\subsection{Discounted Nonnegative Cost Stochastic DP}

\pn We will now consider the  case where $\a<1$.  The cost function of a policy $\p=\{\m_0,\m_1,\ldots\}$ has the form
$$J_\p(x_0)=\lim_{k\to\infty}E_{x_0}^{\p}\lf\{\sum_{ t=0}^{k-1}\a^t g\big(x_ t,\m_ t(x_ t),w_ t\big)\ri\},$$
where as earlier $E^\p_{x_0}\{\cdot\}$ denotes expected value with respect to the probability measure induced by $\p\in \Pi$ under initial state $x_0$. We will assume that  $X$ is a normed space with norm denoted $\|\cdot\|$. 
 
We introduce the set
 $$X_f=\big\{x\in X\mid \jstar(x)<\infty\big\},$$ 
which we assume to be nonempty. Given a state $x\in X_f$, we say that a policy $\p$ is {\it stable from $x$} if there exists a bounded subset of $X_f$ [that depends on $(\p,x)$] such that the (random) sequence $\{x_k\}$ generated starting from  $x$ and using $\pi$ lies with probability 1 within that subset. 
We consider the set
$${\cal C}=\big\{(\p,x)\mid x\in X_f,\ \p\hbox{ is stable from }x\big\},$$
and we assume that ${\cal C}$ is nonempty. 

 Let us say that a function $J\in E^+(X)$ is {\it bounded on bounded subsets of $X_f$\/} if for every bounded subset $\tl X\subset X_f$ there is a scalar $b$ such that $J(x)\le b$ for all $x\in \tl X$. Let us also introduce the set
$$S=\big\{J\in E^+(X)\mid J\hbox{ is bounded on bounded subsets of $X_f$}\big\}.$$
We will assume that $\jstar\in S$. In practical settings we may be able to guarantee this by finding a stationary policy $\m$ such that the function $J_\m$ is bounded on bounded subsets of $X_f$. 
We also assume the following:
 
 \xdef\assumptionstococ{\assumptionn}\assumptionnum\show{myproposition}

\texshopbox{\assumption{\assumptionstococ:}In the discounted stochastic optimal control problem of this section, ${\cal C}$ is nonempty, $\jstar\in S$, and for every $x\in X_f$ and $\e>0$, there exists a policy $\p$ that is stable from $x$ and satisfies $J_\pi(x)\le\jstar(x)+\e$.
}
\smskip

Note that under this assumption, we have $\jstar_{\cal C}=\jstar$. We have the following proposition.

\xdef\propsoc{\propn}\propnum\show{myproposition}

\texshopbox{\proposition{\propsoc:} Let Assumption \assumptionstococ\ hold. Then $\jstar$ is the unique fixed point of $T$ within $S$, and
we have  $T^kJ\to\jstar$ for all $J\in S$ with $\jstar\le J$. If in addition $U$ is a metric space, and the sets $U_k(x,\l)$ of 
Eq.\ \onetwen\ 
are compact for all $x\in X$, $\l\in \re$, and $k$, we have  $T^kJ\to \jstar$ for all $J\in S$, and an optimal stationary policy is guaranteed to exist.}

\proof Using the notation of Section 4.1, we have for all $J\in E(X)$, ${(\p,x_0)}\in{\cal C}$, and $k$,
$$(T_{\m_0}\cdots T_{\m_{k-1}}J)(x_0)=\a^kE_{x_0}^{\p}\big\{J(x_k)\big\}+E_{x_0}^{\p}\lf\{\sum_{ t=0}^{k-1}\a^t g\big(x_ t,\m_ t(x_ t),w_ t\big)\ri\}$$
[cf.\ Eq.\ \finitehorcost]. The fact ${(\p,x_0)}\in{\cal C}$ implies that there is a bounded subset of $X_f$ such that $\{x_k\}$ belongs to that subset with probability 1, so if $J\in S$ it follows  that $\a^kE_{x_0}^{\p}\big\{J(x_k)\big\}\to0$. Thus for all $(\p,x_0)\in{\cal C}$ and $J\in S$, we have
$$\lim_{k\to\infty}(T_{\m_0}\cdots T_{\m_{k-1}}J)(x_0)=\lim_{k\to\infty}E_{x_0}^{\p}\lf\{\sum_{ t=0}^{k-1}\a^t g\big(x_ t,\m_ t(x_ t),w_ t\big)\ri\}=J_\p(x_0),$$
so ${\cal C}$ is $S$-regular. Since $\jstar_{\cal C}$ is equal to $\jstar$ which is a fixed point of $T$ [by Prop.\ \propregset(c)], it follows that $T^kJ\to\jstar$ for all $J\in S$. Under the compactness assumption on the sets $U_k(x,\l)$, the result follows by using Prop.\ \propnegdp(d). \qed

Let us finally note that Assumption \assumptionstococ\ is natural in control contexts where the objective is to keep the state from becoming unbounded, under the influence of random disturbances represented by $w_k$. In such contexts one expects that for a correctly formulated model, optimal or near optimal policies should produce bounded state sequences starting from states with finite optimal cost. 

\vskip-1pc

\section{\bf $S$-Regular Stationary Policies}

\pn We will now specialize the notion of $S$-regularity to stationary policies with the following definition. 

\xdef\definitionregular{\defn}\defnum\show{myproposition}

\texshopbox{\definition{\definitionregular:}For a nonempty set of functions $S\subset E(X)$, we say that a stationary policy $\m$ is  {\it $S$-regular} if  $J_\m\in S$, $J_\m=T_\m J_\m$, and $T^k_\m J\to J_\m$ for all $J\in S$.
A policy that is not $S$-regular is called
{\it $S$-irregular\/}.   
}

Comparing this definition with Definition \definitioncregular, we see that $\m$ is $S$-regular if the set
${\cal C}=\big\{(\m,x)\mid x\in X\big\}$ is $S$-regular, and in addition $J_\m\in S$ and  $J_\m=T_\m J_\m$.
Thus a policy $\m$ is $S$-regular if the VI algorithm corresponding to $\m$,
$J_{k+1}=T_\m J_k,$
represents a dynamic system that has $J_\m$ as its unique equilibrium within $S$, and is asymptotically stable in the sense that the iteration converges to $J_\m$, starting from any $J\in S$.

Generally, with our selection of $S$ we will aim to differentiate between $S$-regular and $S$-irregular policies in a manner that produces useful results for the given problem and does not necessitate restrictive assumptions. 
Examples of sets $S$ that may be fruitfully used are 
$R(X)$, and subsets of $R(X)$ and $E(X)$ involving functions $J$ satisfying $J\ge \jstar$ or $J\ge \bar J$. 
However, there is a diverse range of other useful choices.

\xdef\propregsetcorto{\propn}\propnum\show{myproposition}

\subsection{Restricted Optimization over $S$-Regular Policies}

\pn Given a nonempty set $S\subset E(X)$, let ${\cal M}_S$ be the set of policies that are $S$-regular, and consider optimization over the $S$-regular policies only. The corresponding optimal cost function is denoted $\jstar_S$:
$$\jstar_S(x)=\inf_{\m\in{{\cal M}_S}}J_\m(x),\qquad\forall\ x\in X.\xdef\msinfimum{\lab}\eqnum\show{oneo}$$
We say that $\m^*$ is {\it ${\cal M}_S$-optimal} if 
$$\m^*\in {\cal M}_S\qquad \hbox{and }\qquad J_{\m^*}=\jstar_S.$$
Note that while $S$ is assumed nonempty, it is possible that ${\cal M}_S$ is empty. In this case our results will not be useful, but $\jstar_S$ is still defined by Eq.\ \msinfimum\ as $\jstar_S(x)\equiv\infty$. This is convenient in various proof arguments.

An important question is whether $\jstar_S$ is a fixed point of $T$ and can be obtained by the VI algorithm. Naturally, this depends on the choice of $S$, but it turns out that reasonable choices can be readily found in several important contexts. The following proposition, essentially a specialization of Prop.\ \propregsetcorth,  shows that if  $\jstar_S$ is a fixed point of $T$, then these properties hold within the set
$$W_S=\{J\in E(X)\mid \jstar_S\le J\le \tl J\hbox{ for some $\tl J\in S$}\},\xdef\wellbeh{\lab}\eqnum\show{oneo}$$
which we refer to as the {\it well-behaved region\/}. Note that by the definition of $S$-regularity, the cost functions $J_\m$ of the $S$-regular policies $\m\in {\cal M}_S$ belong to $W_S$. The proposition also provides a necessary and sufficient condition for an $S$-regular policy $\m^*$ to be ${\cal M}_S$-optimal. 

\texshopbox{\proposition{\propregsetcorto:}   Given a set $S\subset E(X)$, assume that $\jstar_S$ is a fixed point of $T$. Then:
\nitem{(a)} ({\it Uniqueness of Fixed Point\/}) $\jstar_S$ is the unique  fixed point of $T$ within $W_S$.%
\nitem{(b)} ({\it VI Convergence\/})  We have $T^kJ\to \jstar_S$ for every
 $J\in W_S$.
\nitem{(c)}  ({\it Optimality Condition\/})  If $\m^*$ is $S$-regular, $\jstar_S \in S$, and $T_{\m^*}\jstar_S=T\jstar_S$, then $\m^*$ is ${\cal M}_S$-optimal. Conversely, if $\m^*$ is ${\cal M}_S$-optimal, then $T_{\m^*}\jstar_S=T\jstar_S$.
}

\proof (a), (b) Follows from Prop.\ \propregsetcorth, with ${\cal C}=\big\{(\m,x)\mid \m\in{\cal M}_S,\,x\in X\big\}$, in which case $\jstar_{\cal C}=\jstar_S$.

\smskip
\pn (c) 
Since $T_{\m^*}\jstar_S=T\jstar_S$ and $T\jstar_S=\jstar_S$, we have $T_{\m^*}\jstar_S=\jstar_S$, and since $\jstar_S\in S$ and $\m^*$ is $S$-regular, we have $\jstar_S=J_{\m^*}$.
Thus $\m^*$ is ${\cal M}_S$-optimal. Conversely, if $\m^*$ is $M_S$-optimal, we have $J_{\m^*}=\jstar_S$, so the fixed point property of $\jstar_S$ and the $S$-regularity of $\m$ imply that 
$T\jstar_S=\jstar_S=J_{\m^*}=T_{\m^*}J_{\m^*}=T_{\m^*}\jstar_S.$ \qed
\vskip-1pc

\xdef \figsimpldecomp{\figr}\figrnum\show{myfigure}

\midinsert
\centerline{\hskip0pc\includegraphics[width=2.5in]{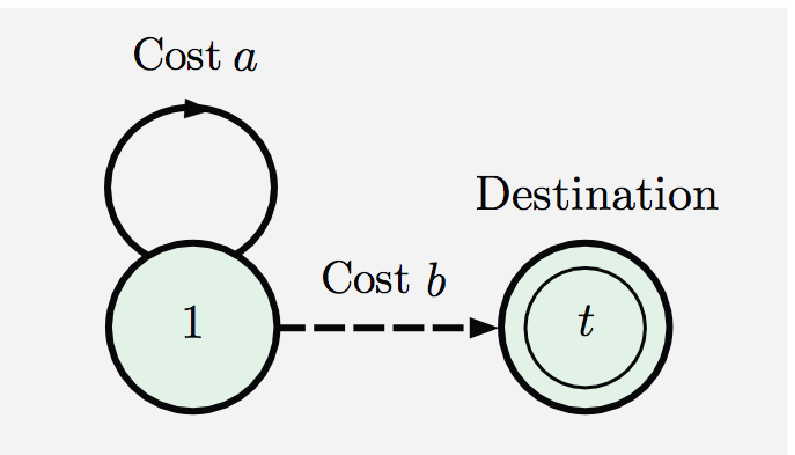}}
\vskip-1pc
\hskip-4pc\fig{0pc}{\figsimpldecomp.} {A shortest path problem with a single node 1 and a termination node $t$.}
\endinsert

The following example illustrates the preceding proposition and demonstrates some of the unusual behaviors that can arise in the context of our model.

\xdef \figdetsp{\figr}\figrnum\show{myfigure}

\xdef\exampledetsp{\exampl}\examplnum\show{myexample}

\beginexample{\exampledetsp}Consider the deterministic shortest path  example shown in Fig.\ \figsimpldecomp. Here there is a single state 1 in addition to the termination state $t$. At state 1 there are two choices: a self-transition, which costs $a$, and a transition to $t$, which costs $b$. The mapping $H$, abbreviating $J(1)$ with just the scalar $J$, is
$$H(1,u,J)=\cases{a+J&if $u$: self transition,\cr
b&if $u$: transition to $t$,\cr}\qquad J\in\re,$$
and the initial function $\skew5\bar J$ is taken to be 0.

There are two policies: the policy $\m$ that transitions from 1 to $t$, which is proper, and the policy $\m'$ that self-transitions at state 1, which is improper. We have
$$T_{\m} J=b,\quad T_{\m'} J=a+J,\quad TJ=\min\{b,\,a+J\},\qquad \forall\ J\in\re.$$
For the proper policy $\m$, the mapping $T_{\m}:\re\mapsto\re$ is a contraction. For the improper policy $\m'$, the mapping $T_{\m'}:\re\mapsto\re$ is not a contraction, and it has a fixed point within $\re$ only if $a=0$, in which case every $J\in\re$ is a fixed point. Let $S$ be equal to the real line $\re$ [the set $R(X)$]. Then a policy is $S$-regular if and only if it is proper (this is generally true for SSP problems,  for $S=\rn$). Thus $\m$ is $S$--regular, while $\m'$ is not. 

Let us consider the optimal cost $J^*$, the fixed points of $T$ within $\re$, and the behavior of VI and PI for different combinations of values of $a$ and $b$.

\nitem{(a)} If  $a>0$, the optimal cost, $J^*=b$, is the unique fixed point of $T$, and the proper policy is optimal.

\nitem{(b)} If $a=0$, the set of fixed points of $T$ (within $\re$) is the interval $(-\infty,b]$. Here the improper policy is optimal if $b\ge0$, and the proper policy is  optimal if $b\le0$. 

\nitem{(c)}  If $a=0$ and $b>0$, the proper policy is strictly suboptimal, yet its cost at state 1 (which is $b$) is a fixed point of $T$. The optimal cost, $J^*=0$, lies in the interior of the set of fixed points of $T$, which is $(-\infty,b]$. Thus the VI method that generates $\{T^kJ\}$ starting with $J\ne J^*$ cannot find $J^*$. In particular if $J$ is a fixed point of $T$, VI stops at $J$, while if $J$ is not a fixed point of $T$ (i.e., $J>b$), VI terminates in two iterations at $b\ne J^*$. Moreover, the standard PI method is unreliable in the sense that starting with the suboptimal proper policy $\m$, it may stop with that policy because $T_{\m}J_{\m}=b=\min\{b,\,J_{\m}\}=TJ_{\m}$ (the improper/optimal policy $\m'$ also satisfies $T_{\m'}J_{\m}=TJ_{\m}$, so a rule for breaking the tie in favor of $\m$ is needed but such a rule may not be obvious in general).

\nitem{(d)} If $a=0$ and $b<0$, the improper policy is strictly suboptimal, and we have $J^*=b$. Here it can be seen that the VI sequence $\{T^kJ\}$ converges to $J^*$ for all $J\ge b$, but stops at $J$ for all $J<b$, since the set of fixed points of $T$ is $(-\infty,b]$. Moreover, starting with either the proper or the improper policy, PI may oscillate, since $T_{\m}J_{\m'}=TJ_{\m'}$ and $T_{\m'}J_{\m}=TJ_{\m}$, as can be easily verified [the optimal policy ${\m}$ also satisfies $T_{\m}J_{\m}=TJ_{\m}$ but it is not clear how to break the tie; compare also with case (c) above]. 

\nitem{(e)} If $a<0$, the improper policy is optimal and we have $J^*=-\infty$. There are no fixed points of $T$ within $\re$, but $J^*$ is the unique fixed point of $T$ within the set $[-\infty,\infty]$.  Then VI will  converge to $J^*$ starting from any $J\in [-\infty,\infty]$, while PI will also converge to the optimal policy starting from either policy.

\smskip
 Let us  focus on the case where there is a zero length cycle ($a=0$).
The cost functions $J_\m$, $J_{\m'}$, and $J^*$ are fixed points of the corresponding mappings, but the sets of fixed points of $T_{\m'}$ and $T$ within $S$ are $\re$ and $(-\infty,b]$, respectively. 
Figure \figdetsp\
 shows the well-behaved regions $W_S$ of Eq.\ \wellbeh\ for the two cases $b>0$ and $b<0$, and is consistent with the results of Prop.\ \propregsetcorto. In particular, the VI algorithm fails when started outside the well-behaved region, while starting from within the region, it is attracted to $J^*_S$ rather than to $J^*$.
\endexample
\vskip-0.5pc

\topinsert
\centerline{\hskip0pc\includegraphics[width=5.5in]{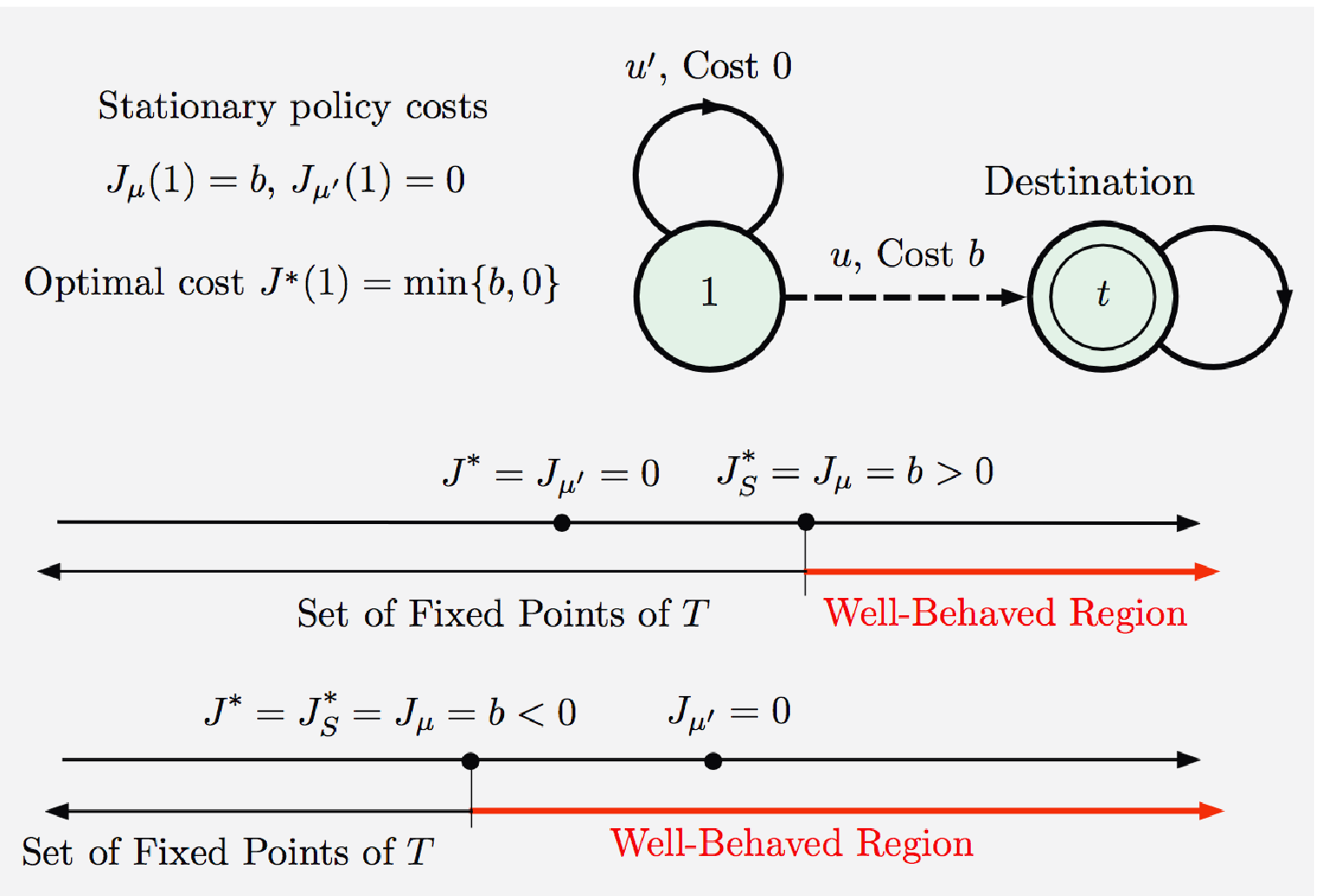}}
\vskip-1pc
\hskip-4pc\fig{0pc}{\figdetsp.} {The well-behaved region of Eq.\ \wellbeh\ for the deterministic shortest path Example \exampledetsp\ when where there is a zero length cycle ($a=0$). For $S=\re$, the policy $\m$ is $S$-regular, while the policy $\m'$ is not. The figure illustrates the two cases where $b>0$ and $b<0$.\vskip-0.5pc}
\endinsert

Note that Prop.\ \propregsetcorto(b) asserts convergence of the VI algorithm to $\jstar_S$ only for initial conditions $J\le \tl J$ for some $\tl J\in S$. For an example where there a single policy $\m$, which is $S$-regular, but $\{T_\m^kJ\}$ does not converge to $J_\m$ starting from some $J\ge J_\m$ that lies outside $S$, consider a mapping $T_\m:\re\mapsto\re$ that has two fixed points: $J_\m$ and another fixed point $J'>J_\m$. Let 
$\tl J=(J_\m+J')/2$
 and  $S=(-\infty,\tl J]$, and assume that $T_\m$ is a contraction mapping within $S$ (a one-dimensional example of this type, where $S=\re$, can be easily constructed graphically). Then, $\tl J\in S$, and starting from any $J\in S$, we have $T^k J\to J_\m$, so that $\m$ is $S$-regular. However, since $J'$ is a fixed point of $T$, the sequence $\{T^kJ'\}$ stays at $J'$ and does not converge to $J_\m$. The difficulty here is that $W_S=[J_\m,\tl J]$ and $J'\notin W_S$.

In many contexts where Prop.\ \propregsetcorto\ applies, there exists an ${\cal M}_S$-optimal policy $\m^*$ such that $T_{\m^*}$ is a contraction with respect to a weighted sup-norm. This is true for example in several types of shortest path problems. In such cases, VI converges to $\jstar_S$ linearly, as shown in the following proposition first given in [BeY16] for SSP problems.
\vskip-0.5pc

\xdef\propvirateconv{\propn}\propnum\show{myproposition}

\texshopbox{\proposition{\propvirateconv: (Convergence Rate of VI)} Let $S$ be equal to $B(X)$, the space of all functions over $X$ that are bounded with respect to a weighted sup-norm $\|\cdot\|_v$ corresponding to a positive function $v:X\mapsto\re$. Assume that $\jstar_S$ is a fixed point of $T$, and that there exists an ${\cal M}_S$-optimal policy $\m^*$ such that  $T_{\m^*}$ is a contraction with respect to $\|\cdot\|_v$,  with corresponding modulus of contraction $\b$. Then
$$\big\|TJ-\jstar_S\|_v\le \b\|J-\jstar_S\|_v,\qquad \forall\ J\ge \jstar_S,\xdef\convrate{\lab}\eqnum\show{oneo}$$
and we have
$$\|J-\jstar_S\|_v\le {1\over 1-\b}\sup_{x\in X}{J(x)-(TJ)(x)\over v(x)},\qquad \forall\ J\ge \jstar_S.\xdef\errbound{\lab}\eqnum\show{oneo}$$
}

\proof By using the {\cal M}-optimality of $\m^*$ and Prop.\ \propregsetcorto(c), we have $\jstar_S=T_{\m^*}\jstar_S=T\jstar_S$, so for all $x\in X$ and $J\ge \jstar_S$,
$${(TJ)(x)-\jstar_S(x)\over v(x)}\le {(T_{\m^*}J)(x)-(T_{\m^*}\jstar_S)(x)\over v(x)}\le \b\max_{x\in X}{J(x)-\jstar_S(x)\over v(x)}.$$
By taking the supremum of the left-hand side over $x\in X$, and by using the fact that the inequality $J\ge \jstar_S$ implies that $TJ\ge T\jstar_S=\jstar_S$, we obtain Eq.\ \convrate.

By using again the relation $T_{\m^*}\jstar_S=T\jstar_S$, we have for all $x\in X$ and all $J\ge \jstar_S$,
$$\eqalign{{J(x)-\jstar_S(x)\over v(x)}&={J(x)-(TJ)(x)\over v(x)}+{(TJ)(x)-\jstar_S(x)\over v(x)}\cr
&\le {J(x)-(TJ)(x)\over v(x)}+{(T_{\m^*}J)(x)-(T_{\m^*}\jstar_S)(x)\over v(x)}\cr
&\le {J(x)-(TJ)(x)\over v(x)}+\b\|J-\jstar_S\|_v.\cr}$$
By taking the supremum of both sides over $x$, we obtain Eq.\ \errbound.
\qed

\subsubsection{Approaches to Show that $\jstar_S$ is a Fixed Point of $T$}

\pn The critical assumption of Prop.\ \propregsetcorto\ is that $\jstar_S$ is a fixed point of $T$. For a specific application, this must be proved with a separate analysis after a suitable set $S$ is chosen. There are several approaches that guide the choice of $S$ and facilitate the analysis. 

One approach applies to problems where $\jstar$ is generically a fixed point of $T$, in which case for every  set $S$ such that $\jstar_S=\jstar$, Prop.\ \propregsetcorto\ applies and shows that $\jstar$ can be obtained by the VI algorithm starting from any $J\in W_S$. This is true generically in wide classes of problems, including deterministic and minimax models (we give a proof for the deterministic case later, in Section 6). 
Other important models where $\jstar$ is guaranteed to be a fixed point of $T$ are the monotone increasing and monotone decreasing models of [Ber13], Section 4.3, a fact known since [Ber77].
In the present paper we will use a different approach for showing that $\jstar_S$ is a fixed point of $T$, which is based on the PI algorithm. 
\vskip-.1pc

\subsection{Policy Iteration-Based Analysis of Bellman's Equation}
\vskip-.3pc
\pn 
In this section we will  develop a PI-based approach for showing that $\jstar_S$ is a fixed point of $T$. The approach is applicable under assumptions that guarantee that there is a sequence $\{\m^k\}$ of $S$-regular policies that can be generated by PI. 
The significance of   all $\m^k$ being $S$-regular  lies in that {\it the corresponding cost function sequence  $\{J_{\m^k}\}$ lies within the well-behaved region of Eq.\ \wellbeh, and is monotonically nonincreasing\/} (see the following Prop.\ 5.2). Under an additional mild technical condition, the limit of this  sequence is a fixed point of $T$ and is in fact equal to $\jstar_S$ (see the subsequent Prop.\ 5.3).

Let us consider the PI algorithm that generates a sequence of policies $\{\m^k\}$ according to
$$T_{\m^{k+1}}J_{\m^k}=T J_{\m^k},\qquad k=0,1,\ldots,\xdef\poiter{\lab}\eqnum\show{oneo}$$
starting from an initial policy $\m^0$. This iteration embodies both the policy evaluation step, which computes $J_{\m^k}$ in some way, and the policy improvement step, which computes $\m^{k+1}(x)$ as a minimum over $u\in U(x)$ of $H(x,u,J_{\m^k})$ for each $x\in X$ [cf.\ Eq.\ \poiter]. Of course,  to be able to carry out the policy improvement step, there should be enough assumptions to guarantee that the minimum is attained for every $x$. One such assumption is that $U(x)$ is a finite set for each $x\in X$. A more general assumption, involving a form of compactness of the constraint set is given in the next section (see Lemma 6.1).

The  evaluation of the cost function $J_\m$ of a policy $\m$ may be done by solving the equation $J_\m=T_\m J_\m$, which holds when $\m$ is an $S$-regular policy. An important fact is that if the PI algorithm generates a sequence $\{\m^k\}$ consisting exclusively of $S$-regular policies, then not only the policy evaluation is facilitated through the equation $J_\m=T_\m J_\m$, but also {\it the sequence of  cost functions $\{J_{\m^k}\}$ is monotonically nonincreasing\/}, as we will show next. 

Note a fine point here. For a given starting policy $\m^0$, there may be many different sequences $\{\m^k\}$ that can be generated by  PI [i.e., satisfy Eq.\ \poiter]. Some of these may consist  of $S$-regular policies exclusively, and some may not. The policy improvement property shown in the following proposition holds for the former sequences, but not necessarily for the latter.

\xdef\propimprovereg{\propn}\propnum\show{myproposition}

\texshopbox{\proposition{\propimprovereg: (Policy Improvement Under $S$-Regularity)}Given a set $S\subset E(X)$, assume that $\{\m^k\}$ is a sequence  generated by the PI algorithm \poiter\ that consists of $S$-regular policies. Then
$J_{\m^k}\ge J_{\m^{k+1}}$ for all  $k$.
}

\proof Using the $S$-regularity of $\m^k$, we have
$$J_{\m^k}=T_{\m^k}
J_{\m^k}\ge TJ_{\m^k}= T_{\m^{k+1}}J_{\m^k}.\xdef\costimproveo{\lab}\eqnum\show{oneo}$$ 
By repeatedly applying $T_{\m^{k+1}}$ to both sides, we obtain
$$J_{\m^k}\ge \lim_{m\tends\infty}T_{\m^{k+1}}^m J_{\m^k}=J_{\m^{k+1}},$$
where the equation on the right holds since ${\m^{k+1}}$ is $S$-regular and $J_{\m^k}\in S$ (since $\m^k$ is $S$-regular). \qed

The preceding proposition shows that if a sequence of $S$-regular policies $\{{\m^k}\}$ is generated by PI, the corresponding  cost function sequence $\{J_{\m^k}\}$ is monotonically nonincreasing and hence converges to a limit $J_\infty$. Under mild conditions, we will show that $J_\infty$ is a fixed point of $T$ and is equal to $\jstar_S$. This is important as it brings to bear Prop.\ \propregsetcorto, and the associated results on VI convergence and optimality conditions. 
Let us first formalize the property that the PI algorithm can generate a sequence of $S$-regular policies.

\xdef\definitionpiclosureweak{\defn}\defnum\show{myproposition}

\texshopbox{\definition{\definitionpiclosureweak: (Weak PI  Property)}A set $S\subset E(X)$ has the {\it weak PI property} if there exists a sequence of $S$-regular policies that can be generated by the PI algorithm [i.e., a sequence $\{\m^k\}$ that satisfies Eq.\ \poiter\ and consists of $S$-regular policies]. 
}

The following proposition provides the basis for showing that $\jstar_S$ is a fixed point of $T$ based on the weak PI  property and a mild continuity-type condition. 

\xdef\propdetpolicyitgen{\propn}\propnum\show{myproposition}

\texshopbox{\proposition{\propdetpolicyitgen: (Weak PI Property Theorem)}Given a set $S\subset E(X)$, assume that:
\nitem{(1)} $S$ has the weak PI property.
\nitem{(2)} For each sequence  $\{J_m\}\subset  S$ with $J_m\downarrow J$ for some $J\in E(X)$, we have 
$$H\lf(x,u,J\ri)= \lim_{m\to\infty}H(x,u,J_m),\qquad \forall\ x\in X,\ u\in U(x).\xdef\limrelation{\lab}\eqnum\show{oneo}$$
\pn 
Then:
\nitem{(a)}  $\jstar_S$ is a fixed point of $T$ and the  conclusions of Prop.\ \propregsetcorto\ hold. 
\nitem{(b)} ({\it PI Convergence\/}) Every sequence of $S$-regular policies $\{\m^k\}$ that can be generated by PI satisfies $J_{\m^k}\downarrow \jstar_S$.  If in addition the set of $S$-regular policies is finite, there exists $\bar k\ge 0$ such that $\m^{\bar k}$ is ${\cal M}_S$-optimal.
}

\proof (a) Let $\{\m^k\}$ be a sequence  of $S$-regular policies generated by the PI algorithm (there exists such a sequence by the weak PI property). Then by Prop.\ \propimprovereg, the sequence $\{J_{\m^k}\}$ is monotonically nonincreasing and must converge to some $J_\infty\ge \jstar_S$. 
We will show that $J_\infty$ is a fixed point of $T$ and then invoke Prop.\ \propregsetcorth.

Indeed, we have 
$$J_{\m^k}\ge TJ_{\m^k}\ge TJ_\infty$$
 [cf.\ Eq.\ \costimproveo], so by letting $k\to\infty$, we obtain $J_\infty\ge TJ_\infty$.
To prove the reverse inequality, we first note that from the definition of the PI iteration and the  nonincreasing property $J_{\m^k}\ge J_{\m^{k+1}}$, we have
$$T J_{\m^k}=T_{\m^{k+1}}J_{\m^k}\ge T_{\m^{k+1}}J_{\m^{k+1}}=J_{\m^{k+1}}.$$
By using Eq.\ \limrelation\ together with the preceding relation, we obtain for all $x\in X$ and $u\in U(x)$,
$$H(x,u,J_\infty)= \lim_{k\to\infty}H(x,u,J_{\m^k})\ge \lim_{k\to\infty}\,(TJ_{\m^k})(x)\ge \lim_{k\to\infty}\,J_{\m^{k+1}}=J_\infty(x).$$
By taking the infimum of the left-hand side over $u\in U(x)$, it follows that $TJ_\infty\ge J_\infty$.  Thus $J_\infty=TJ_\infty$.
Finally, by applying Prop.\ \propregsetcorth\ with ${\cal C}=\big\{(\m,x)\mid \m\in{\cal M}_S,\,x\in X\big\}$, we have $J_\infty=\jstar_{\cal C}=\jstar_S$.
\smskip
\pn {(b)} The limit of $\{J_{\m^k}\}$ was shown to be equal to $\jstar_S$ in the preceding proof. Moreover, the finiteness of ${\cal M}_S$ and  the policy improvement property of Prop.\ \propimprovereg\ imply that some $\m^{\bar k}$ is ${\cal M}_S$-optimal.
\qed

Note that under the weak PI property, the preceding proposition shows convergence of PI to $\jstar_S$ but not necessarily to $\jstar$. Moreover, it is possible for the PI algorithm to generate a nonmonotonic sequence of policy cost functions that includes both optimal and strictly suboptimal policies, as was seen in the deterministic shortest path Example \exampledetsp\ for the case where $a=0$ and $b<0$.

Proposition \propdetpolicyitgen(a) does not guarantee that {\it every} sequence $\{\m^k\}$ generated by the PI algorithm satisfies $J_{\m^k}\downarrow \jstar_S$. This is true only for the sequences that consist of $S$-regular policies. We know that when the weak PI property holds, there exists at least one such sequence, but PI can also generate sequences that contain $S$-irregular policies, as we have seen in Example \exampledetsp. We thus introduce a stronger type of PI property, which we will use to obtain stronger results.

\xdef\definitionpiclosurestrong{\defn}\defnum\show{myproposition}

\texshopbox{\definition{\definitionpiclosurestrong: (Strong PI  Property)}A set $S\subset E(X)$ has the {\it strong PI property} if:
\nitem{(a)} There exists at least one $S$-regular policy.
\nitem{(b)} For every $S$-regular policy $\m$, any policy $\bar \m$ such that $T_{\bar \m}J_\m=TJ_\m$ is $S$-regular, and there exists at least one such $\bar \m$.
}

The strong PI property implies that every sequence that can be generated by PI starting from an $S$-regular policy consists exclusively of $S$-regular policies. Moreover, there exists at least one such sequence. Hence the strong PI property implies the weak PI property. Thus if the strong PI property holds together with the mild continuity condition (2) of Prop.\ \propdetpolicyitgen, $\jstar_S$ is a fixed point of $T$ and Prop.\ \propregsetcorto\ applies. 

On the other hand, the strong PI property may be harder to verify in a given setting. The following proposition provides conditions guaranteeing that $S$ has the strong PI property. The key implication of these conditions is that they preclude optimality of an $S$-irregular policy [see condition (4) of the proposition]. Condition (3) of the proposition is implied by finiteness of the constraint set or by a more general compactness assumption that will be given in the next section.

\xdef\propdetpolicyitgenth{\propn}\propnum\show{myproposition}

\texshopboxnb{\proposition{\propdetpolicyitgenth: (Verifying the Strong PI Property)}Given a set $S\subset E(X)$, assume that:
\nitem{(1)} $J(x)<\infty$ for all $J\in S$ and $x\in X$.
\nitem{(2)} There exists at least one $S$-regular policy.
\nitem{(3)} For every $J\in S$ there exists a policy $\m$ such that $T_{\m}J=TJ$.
\nitem{(4)} For every $J\in S$ and $S$-irregular policy $\m'$, there exists  a state $x\in X$ such that}\texshopboxnt{\nitem{}
$$\limsup_{k\to\infty}\,(T_{\m'}^k J)(x)=\infty.\xdef\infcondition{\lab}\eqnum\show{oneo}$$
\pn Then:
\nitem{(a)} If a policy $\m$ satisfies $T_\m J\le J$ for some function $J\in S$, then $\m$ is $S$-regular. 
\nitem{(b)} $S$ has the strong PI  property.
}

\proof (a) By the monotonicity of $T_\m$, we have $\limsup_{k\to\infty}T_\m^k J\le J$, and since by condition (1), $J(x)<\infty$ for all $x$, it follows from Eq.\ \infcondition\ that $\m$ is $S$-regular. 
\smskip
\pn (b) In view of condition (3), it will suffice to show that for every $S$-regular policy $\m$, any policy $\bar \m$ such that $T_{\bar \m}J_\m=TJ_\m$ is also $S$-regular. Indeed we have
$T_{\bar \m}J_\m=TJ_\m\le T_\m J_\m=J_\m,$
so $\bar \m$ is $S$-regular by part (a). \qed

By using the strong PI property and assuming also that $\jstar_S\in S$, we will now show that $\jstar_S$ is the unique fixed point of $T$ within $S$. This result will be the starting point for the analysis of Section 6.

\xdef\propdetpolicyitgenf{\propn}\propnum\show{myproposition}

\texshopbox{\proposition{\propdetpolicyitgenf: (Strong PI Property Theorem)}Let $S$ satisfy the conditions of Prop.\ \propdetpolicyitgenth.
\nitem{(a)} ({\it Uniqueness of Fixed Point\/}) If $T$ has a fixed point within $S$, then this fixed point is equal to $\jstar_S$. 
\nitem{(b)} ({\it Fixed Point Property and Optimality Condition\/}) If $\jstar_S\in S$, then $\jstar_S$ is the unique fixed point of $T$ within $S$. Moreover, every policy $\m$ that satisfies $T_\m \jstar_S=T\jstar_S$ is ${\cal M}_S$-optimal and there exists at least one such policy.
\nitem{(c)}  ({\it PI Convergence\/}) If for each sequence  $\{J_m\}\subset  S$ with $J_m\downarrow J$ for some $J\in E(X)$, we have 
$$H\lf(x,u,J\ri)= \lim_{m\to\infty}H(x,u,J_m),\qquad \forall\ x\in X,\ u\in U(x),\eqnum\show{oneo}$$
then  $\jstar_S$ is a fixed point of $T$, and every sequence $\{\m^k\}$ generated by the PI algorithm starting from an $S$-regular policy $\m^0$ satisfies $J_{\m^k}\downarrow\jstar_S$. Moreover, if the set of $S$-regular policies is finite, there exists $\bar k\ge 0$ such that $\m^{\bar k}$  is ${\cal M}_S$-optimal.}

\proof (a) Let $J'\in S$ be a fixed point of $T$. By applying Prop.\ \propregsetcorth\ with ${\cal C}=\big\{(\m,x)\mid \m\in{\cal M}_S,\,x\in X\big\}$, we have $J'\le \jstar_{\cal C}=\jstar_S$. 
For the reverse inequality, let $\m'$ be such that $J'=TJ' =T_{\m'}J'$ [cf.\ condition (3) of Prop.\ \propdetpolicyitgenth]. Then by Prop.\ \propdetpolicyitgenth(a), it follows that $\m'$ is $S$-regular, and since $J'\in S$, by the definition of $S$-regularity, we have 
$J' =J_{\m'}\ge \jstar_S$, showing that $J'=\jstar_S$. 
\smskip

\pn (b) For every $\m\in{\cal M}_S$ we have $J_\m\ge \jstar_S$, so that
$J_\m=T_\m J_\m\ge T_\m \jstar_S\ge T \jstar_S.$
Taking the infimum over all $\m\in{\cal M}_S$, we obtain
$\jstar_S\ge T \jstar_S.$
Let $\m$ be a policy such that 
$T \jstar_S=T_\m\jstar_S,$
[there exists one by condition (3) of Prop.\ \propdetpolicyitgenth, since we assume that $\jstar_S\in S$]. The preceding two relations yield
$\jstar_S\ge T_\m \jstar_S$, so by Prop.\ \propdetpolicyitgenth(a), $\m$ is $S$-regular. Therefore, we have
$$\jstar_S\ge T \jstar_S= T_\m \jstar_S\ge \lim_{k\to\infty}T_\m^k \jstar_S=J_\m\ge \jstar_S,$$
where the second equality holds by $S$-regularity of $\m$ and $\jstar_S\in S$ by assumption. Hence equality holds throughout in the above relation, proving that $\jstar_S$ is a fixed point of $T$ and that $\m$ is ${\cal M}_S$-optimal.

\smskip
\pn (c) Since the strong PI property [which holds by Prop.\ \propdetpolicyitgenth(b)] implies the weak PI property, the result follows from Prop.\ \propdetpolicyitgen. 
\qed

The preceding proposition does not address the question whether $\jstar$ is a fixed point of $T$, and does not guarantee that VI converges to $\jstar_S$ or $\jstar$ starting from every $J\in S$. We will consider both of these issues in the next section.  Note a simple consequence of part (a): if $\jstar$ is known to be a fixed point of $T$ and to belong to $S$, then $\jstar=\jstar_S$. 
\old{
\xdef\exampledetspthree{\exampl}\examplnum\show{myexample}
\beginexample{\exampledetspthree: (Strong PI Property and the Deterministic Shortest Path Example)}Consider the deterministic shortest path Example \exampledetsp\ for the case where the cycle has positive length ($a>0$), and let $S$ be the real line $\re$. The two policies are: $\m$ which moves from state 1 to the destination at cost $b$ and is $S$-regular, and $\m'$ which stays at state 1 at cost $a$, which is $S$-irregular. However, $\m'$ has infinite cost and satisfies Eq\ \infcondition. As a result, Prop.\ \propdetpolicyitgenth\ applies and the strong PI property holds. Consistent with Prop.\ \propdetpolicyitgenf, $J^*_S$ is the unique fixed point of $T$ within $S$. 
Turning now to  the PI algorithm, we see that starting from the $S$-regular $\m$, which is optimal,  it stops at $\m$, consistent with Prop.\ \propdetpolicyitgenf(e). However, starting from the $S$-irregular policy $\m'$ the policy evaluation portion of the PI algorithm must be able to deal with the infinite cost values associated with $\m'$. This is a generic difficulty in applying PI to problems where there are irregular policies: we either need to know an initial $S$-regular policy, or appropriately modify the PI algorithm.
\endexample
}

Proposition \propdetpolicyitgenf(c) shows that PI is valid, but for this an initial $S$-regular policy must be available. 
Chapter 3 of [Ber13] describe a combined VI and PI algorithm, which does not require an initial $S$-regular policy, and can tolerate the generation of $S$-irregular policies. Let us also consider two additional algorithmic approaches for computing $\jstar_S$, not given in [Ber13], which can be justified based on the preceding analysis.

\subsubsection{A Mathematical Programming Solution Method}

\pn We will  show that $\jstar_S$ is an upper bound to all functions $J\in S$ that satisfy $J\le TJ$, and we will exploit this fact to obtain an algorithm to compute $\jstar_S$. We have the following proposition. 

\xdef\proplinprogr{\propn}\propnum\show{myproposition}

\texshopbox{\proposition{\proplinprogr:}  Given a set $S\subset E(X)$,  for all functions $J\in S$ satisfying $J\le TJ$, we have $J\le \jstar_S$.}

\proof If $J\in S$ and $J\le TJ$, by repeatedly applying $T$ to both sides and using the monotonicity of $T$, we obtain
$J\le T^kJ\le T_\m^k J$ for all $k$ and $S$-regular policies $\m$.
Taking the limit as $k\to\infty$, we obtain $J\le J_\m$, so by taking the infimum over $\m\in{\cal M}_S$, we obtain $J\le \jstar_S$. \qed

Assuming that $\jstar_S$ is a fixed point of $T$, we can use the preceding proposition to compute $\jstar_S$ by maximizing an appropriate monotonically increasing function of $J$ subject to the constraints $J\in S$ and $J\le TJ$.\footnote{$\,$\dag}{\ninepoint For the mathematical programming approach to apply, it is sufficient that $J^*_S\le T J^*_S$. However, we generally have $J^*_S\ge TJ^*_S$ (this follows by writing for all $\m\in {\cal M}_S$, $J_\m=T_\m J_\m\ge TJ_\m\ge TJ^*_S$, and taking the infimum over all $\m\in {\cal M}_S$), so the condition $J^*_S\le TJ^*_S$ is equivalent to $J^*_S$ being a fixed point of $T$.} This approach is well-known in finite-state finite-control Markovian decision problems, where it is usually referred to as the {\it linear programming solution method\/}, because in this case the resulting optimization problem is a linear program  (see e.g., the books [Kal83], [Put94], [Ber12]).

 For a more general finite-state case, suppose that $X=\{1,\ldots,n\}$ and $S=\rn$. Then Prop.\ \proplinprogr\ shows that $\jstar_S=\big(\jstar_S(1),\ldots,\jstar_S(n)\big)$ is the unique solution of the following optimization problem:
$$\eqalign{\hbox{\rm maximize}\quad &\sum_{i=1}^n \b_i J(i)\cr
\hbox{\rm subject to\ \ }
&J(i)\le H(i,u,J),\ \ i =1,\ldots,n, \quad u\in U(i),\cr}\old{\eqnum\show{oneo}}$$
where $\b_1,\ldots,\b_n$ are any positive scalars.
If $H$ is linear in $J$ and each $U(i)$ is a finite set, this is a linear program, which can be solved by using standard linear programming methods.

\subsubsection{An Optimistic Form of PI}

\pn Let us finally consider an optimistic variant of PI, where policies are evaluated inexactly, with a finite number of VIs. In particular,
this algorithm starts with some $J_0\in E(X)$ such that $J_0\ge TJ_0$, and generates a sequence $\{J_k,\m^k\}$ according to
$$T_{\m^k} J_k=T J_k,\qquad J_{k+1}=T_{\m^k}^{m_k}J_k,\qquad k=0,1,\ldots,\xdef\optpoiter{\lab}\eqnum\show{oneo}$$
where $m_k$ is a positive integer for each $k$. 

The following proposition shows that optimistic PI converges under mild assumptions to a fixed point of $T$, independently of any $S$-regularity framework. However, when such a framework is introduced, and the sequence generated by optimistic PI generates a sequence of $S$-regular policies, then the algorithm converges to $\jstar_S$, which is in turn a fixed point of $T$, similar to the PI convergence result under the weak PI property; cf.\ Prop.\ \propdetpolicyitgen(b). Thus the proposition serves both an analytical purpose (as a tool for establishing that $\jstar_S$ is a fixed point of $T$), and a computational purpose [establishing the validity of the optimistic PI algorithm \optpoiter\ as a means for computing $\jstar_S$].

\xdef\propdetpolicyitopt{\propn}\propnum\show{myproposition}

\texshopbox{\proposition{\propdetpolicyitopt: (Convergence of Optimistic PI)} Let $J_0\in E(X)$ be a function such that $J_0\ge TJ_0$, and assume that:
\nitem{(1)} For all $\m\in{\cal M}$, we have $J_\m=T_\m J_\m$,  and for all $J\in E(X)$ with $J\le J_0$, there exists $\bar\m\in{\cal M}$ such that $T_{\bar\m}J=TJ$.
\nitem{(2)} For each sequence  $\{J_m\}\subset  E(X)$ with $J_m\downarrow J$ for some $J\in E(X)$, we have
$$H\lf(x,u,J\ri)= \lim_{m\to\infty}H(x,u,J_m),\quad \forall\ x\in X,\ u\in U(x).$$
\pn Then the optimistic PI algorithm \optpoiter\ is well defined and  the following hold:
\nitem{(a)} The sequence $\{J_k\}$ generated by the algorithm satisfies $J_{k}\downarrow J_\infty$, where $J_\infty$ is a fixed point of $T$.
\nitem{(b)} If for a set $S\subset E(X)$, the sequence $\{\m^k\}$ generated by the algorithm consists of $S$-regular policies and we have $J_k\in S$ for all $k$, then $J_k\downarrow \jstar_S$ and $\jstar_S$ is a fixed point of $T$.
}

 \proof (a) Condition (1) guarantees that the  sequence $\{J_k,\m^k\}$ is well defined in the following argument. We also have
$$J_0\ge TJ_0= T_{\m^0}J_0\ge T_{\m^0}^{m_0}J_0= J_1\ge T_{\m^0}^{m_0+1}J_0=T_{\m^0}J_1\ge TJ_1=T_{\m^1}J_1\ge\cdots\ge J_2,\xdef\costimproveopt{\lab}\eqnum\show{oneo}$$
and continuing similarly, we obtain 
$J_k\ge TJ_k\ge J_{k+1}$ for all  $k=0,1,\ldots.$
Thus  $J_{k}\downarrow J_\infty$ for some $J_\infty$. The proof that $J_\infty$ is a fixed point of $T$ is the same as in the case of the PI algorithm \poiter\ in Prop.\ \propdetpolicyitgen. 
\smskip 
\pn (b) In the case where all the policies $\m^k$ are $S$-regular and $\{J_k\}\subset S$, from Eq.\ \costimproveopt, we have $J_{k+1}\ge J_{\m^k}$ for all $k$, so it follows that
$$J_\infty=\lim_{k\to\infty}J_k\ge \liminf_{k\to\infty}J_{\m^k}\ge \jstar_S.$$
We will also show that the reverse inequality holds, so that $J_\infty=\jstar_S$. Indeed, for every $S$-regular policy $\m$ and all $k\ge0$, we have 
$$J_\infty=T^kJ_\infty\le T_\m^k J_\infty\le  T_\m^k J_0,$$
from which by taking limit as $k\to\infty$ and using the assumption $J_0\in S$, we obtain 
$$J_\infty\le \lim_{k\to\infty}T_\m^k J_0=J_\m,\qquad \forall\ \m\in {\cal M}_S.$$
Taking the infimum over $\m\in {\cal M}_S$, it follows that $J_\infty\le \jstar_S$. 
Thus, $J_\infty=\jstar_S$, and by using the properties of $J_\infty$ proved in part (a), the result follows. \qed

Note that the fixed point $J_\infty$ in Prop.\ \propdetpolicyitopt(a) need not be equal to $\jstar_S$ or $\jstar$. As an illustration, consider the shortest path Example \exampledetsp\ with $S=\re$, and $a=0$, $b>0$. Then if $0<J_0<b$, it can be seen that $J_k=J_0$ for all $k$, so $\jstar=0<J_\infty$ and  $J_\infty<\jstar_S=b$.

\vskip-0.5pc\vskip-0.5pc
\section{Irregular Policies/Infinite Cost Case}
\vskip-0.5pc

\pn The results of the preceding section do not assert that $\jstar$ is a fixed point of $T$ or that $\jstar=\jstar_S$. In this section we address this issue with some additional assumptions. The following assumption and proposition were first given in Section 3.2 of [Ber13], but the line of proof given here is considerably streamlined thanks to the use of the strong PI property analysis of the preceding section, which was developed after [Ber13] was published.
\vskip-0.5pc
\xdef\assumptiontoz{\assumptionn}\assumptionnum\show{myproposition}

\texshopboxnb{
\assumption{\assumptiontoz:} We have a subset $S\subset R(X)$ satisfying the following:
\nitem{(a)} $S$ contains $\bar J$, and has the property that if $J_1,J_2$ are two functions in $S$, then $S$ contains all functions $J$ with $J_1\le J\le J_2$.%
\nitem{(b)} The function $\jstar_S=\inf_{\m\in{\cal M}_S}J_\m$
belongs to $S$.
\nitem{(c)} For each $S$-irregular policy $\m$ and each $J\in S$, there is at least one state $x\in X$ such that  
$$\limsup_{k\to\infty}\,(T_\m^kJ)(x)=\infty.\eqnum\show{oneo}$$}\texshopboxnt{\nitem{}
\nitem{(d)} The control set $U$ is a metric space, and the set
$$
\big\{ u\in U(x)\mid H(x,u,J)\le
\l\big\}$$
is compact for every
$J\in S$, $x\in X$, and $\l\in \re$.
\nitem{(e)} For each sequence  $\{J_m\}\subset S$ with $J_m\uparrow J$ for some $J\in S$,
$$\lim_{m\to\infty}H(x,u,J_m)= H\lf(x,u,J\ri),\qquad \forall\ x\in X,\ u\in U(x).$$
\nitem{(f)} For each  function $J\in S$, there exists a function  $J'\in S$ such that $J'\le J$ and $J'\le TJ'$.
}
\vskip-0.2 pc

The  conditions (b) and (c) of the preceding assumption have been introduced in Props.\ \propdetpolicyitgenth\ and \propdetpolicyitgenf\ in the context of the strong PI property-related analysis. New conditions, not encountered earlier, are (a), (d), (e), and (f). They will be used to  assert that $\jstar=\jstar_S$, that $\jstar$ is the unique fixed point of $T$ within $S$, and that the VI and PI algorithms have improved convergence properties compared with the ones of Section 5.2, thereby obtaining results that are almost as strong as the ones of Chapter 2 for contractive models. In the case where $S$ is the set of real-valued functions $R(X)$ and $\bar J\in R(X)$, condition (a) is automatically satisfied, while condition (e) is typically verified easily. The verification of condition (f) may be nontrivial in some cases. We postpone the discussion of this issue for later (see the subsequent Prop.\ 6.2).

The main result of this section is the following proposition. 

\xdef\propttffalter{\propn}\propnum\show{myproposition}

\texshopbox{\proposition{\propttffalter:}Let Assumption \assumptiontoz\ hold.
Then:
\nitem{(a)} The optimal cost function $\jstar$ is the unique fixed point of $T$ within the set $S$.
\nitem{(b)} We have $T^k J\to \jstar$ for all $J\in S$.
\nitem{(c)} A policy  $\m$ is optimal if and only if
$T_\m \jstar =T\jstar$.  Moreover, there exists an optimal $S$-regular policy.
\nitem{(d)} For any $J\in S$, if $J\le TJ$ we have $J\le \jstar $, and if $J\ge TJ$ we have $J\ge \jstar $.
\nitem{(e)} If in addition for each sequence  $\{J_m\}\subset  S$ with $J_m\downarrow J$ for some $J\in S$, we have 
$$H\lf(x,u,J\ri)= \lim_{m\to\infty}H(x,u,J_m),\qquad \forall\ x\in X,\ u\in U(x),\eqnum\show{oneo}$$
then every sequence $\{\m^k\}$ generated by the PI algorithm starting from an $S$-regular policy $\m^0$ satisfies $J_{\m^k}\downarrow\jstar$. Moreover, if the set of $S$-regular policies is finite, there exists $\bar k\ge 0$ such that $\m^{\bar k}$ is optimal.
}

The proof of Prop.\ \propttffalter\ will make use of the analysis of the preceding section. We first state without proof a result given as Lemma 3.2.1 of [Ber13]. It guarantees that  starting from an $S$-regular policy, the PI algorithm is well defined. Similar results are well-known in DP theory.

\xdef\lemmao{\lemman}\lemmanum\show{myproposition}

\texshopbox{
\lemma{\lemmao:} Let Assumption \assumptiontoz(d) hold. For every $J\in S$, there exists a policy $\m$ such that $T_\m J=TJ$.
}
\vskip-0.5pc

Next we restate, for easy reference,  some of the results of the preceding section in the next two lemmas. 

\xdef\lemmat{\lemman}\lemmanum\show{myproposition}

\texshopbox{
\lemma{\lemmat:} Let Assumption \assumptiontoz(c) hold. A policy $\m$ that satisfies $T_\m J\le J$ for some $J\in S$ is $S$-regular.
}

\proof This is Prop.\ \propdetpolicyitgenth(b).  \qed
\vskip-0.5pc
\vskip-0.5pc

\xdef\lemmath{\lemman}\lemmanum\show{myproposition}

\texshopbox{
\lemma{\lemmath:} Let Assumption \assumptiontoz(b),(c),(d) hold.
Then:
\nitem{(a)} The function $\jstar_S$ of Assumption \assumptiontoz(b)
is the unique fixed point of $T$ within $S$.
\nitem{(b)} Every policy $\m$ satisfying $T_{\m}\jstar_S=T\jstar_S$ is optimal within the set of $S$-regular policies, i.e., $\m$ is $S$-regular and $J_{\m}=\jstar_S$. Moreover, there exists at least one such policy.
}

\proof This is Prop.\ \propdetpolicyitgenf, parts (a) and (b)  [Assumption \assumptiontoz(d) guarantees that for every $J\in S$, there exists a policy $\m$ such that $T_\m J=TJ$ (cf.\ Lemma \lemmao)].
\qed

Let us also prove the following technical lemma that relies on the continuity Assumption \assumptiontoz(e).

\xdef\lemmaf{\lemman}\lemmanum\show{myproposition}

\texshopbox{
\lemma{\lemmaf:} Let Assumption \assumptiontoz(d),(e) hold. Then if $J\in S$, $\{T^kJ\}\subset S$, and $T^kJ\uparrow J_\infty$ for some $J_\infty\in S$, we have $J_\infty= \jstar_S$.}

\proof We fix $x\in X$, and consider the sets 
$$U_k(x)
= \Big\{ u\in U(x)
\mid H(x,u,T^kJ)\le J_\infty(x)\Big\},\qquad k=0,1,\ldots,\xdef\ukdef{\lab}\eqnum\show{oneo}$$ 
which are compact by assumption. Let
 $u_k\in U(x)$ be such that  
$$H(x,u_k,T^kJ) = \inf_{u\in U(x)}H(x,u,T^kJ)=(T^{k+1}J)(x)\le J(x)$$ 
(such a point exists by Lemma \lemmao). Then $u_k\in U_k(x)$.

For every
$k$, consider the sequence $\{ u_i\}^\infty_{i=k}$.  Since 
$T^kJ\uparrow J_\infty,$
 it follows that for all $i\ge k$,
$$H(x,u_i,T^kJ)\le H(x,u_i,T^iJ)\le J_\infty(x).$$ 
Therefore from the definition \ukdef, we have $\{
u_i\}_{i=k}^\infty \subset U_k(x)$. Since $U_k(x)$ is compact, all the limit points of $\{ u_i\}_{i=k}^\infty$ belong to
$U_k(x)$ and at least one limit point exists.  Hence the
same is true for the limit points of the whole sequence $\{ u_i\}$.  Thus if $\tilde u$ is a limit point of $\{ u_i\}$, we have 
$$\tilde
u\in \cap_{k=0}^\infty U_k(x).$$ 
By Eq.\ \ukdef, this implies that 
$$H\big(x,\tilde u,T^kJ\big)\le J_\infty(x),\qquad k=0,1,\ldots.$$
Taking
the limit as $k\tends \infty$ and using Assumption \assumptiontoz(e), we obtain 
$$(T J_\infty)(x)\le H(x,\tilde u,J_\infty)\le J_\infty(x).$$
Thus, since $x$ was chosen arbitrarily within $X$, we have $TJ_\infty\le J_\infty$. To show the reverse inequality, we write $T^kJ\le J_\infty$, apply $T$ to this inequality, and take the limit as $k\to\infty$, so that $J_\infty=\lim_{k\to\infty}T^{k+1}J\le TJ_\infty$. It follows that  $J_\infty=TJ_\infty$. Since $J_\infty\in S$, by part (a) we have $J_\infty=\jstar_S$. \qed

We are now ready to show Prop.\ \propttffalter\ by using the additional parts (a) and (f) of Assumption \assumptiontoz.

\smskip

\pn {\bf Proof of  Prop.\ \propttffalter:} (a), (b) We will first prove that $T^kJ\tends \jstar_S$ for all $J\in S$, and we will use this to prove  that $\jstar_S=\jstar$ and that there exists an optimal $S$-regular policy. Thus  parts (a)  and (b), together with the existence of an optimal $S$-regular policy, will be shown simultaneously.
 
We fix $J\in S$, and choose $J'\in S$ such that $J'\le J$ and $J'\le TJ'$ [cf.\ Assumption \assumptiontoz(f)]. By the monotonicity of $T$, we have $T^kJ'\uparrow J_\infty$ for some $J_\infty\in E(X)$. Let $\m$ be an $S$-regular policy such that $J_{\m}=\jstar_S$ [cf.\ Lemma \lemmath(b)]. Then we have, using again  the monotonicity of $T$,
$$J_\infty=\lim_{k \to\infty}T^k J'\le \limsup_{k\to\infty}T^kJ\le \lim_{k\to\infty}T_\m^k J=J_\m=\jstar_S.\xdef\jtildejhat{\lab}\eqnum\show{oneo}$$
Since $J'$ and $\jstar_S$ belong to $S$, and $J'\le T^kJ'\le  J_\infty\le \jstar_S$, Assumption \assumptiontoz(a) implies that $\{T^kJ'\}\subset S$, and $J_\infty\in S$. From Lemma \lemmaf, it then follows that $J_\infty=\jstar_S$.  Thus equality holds throughout in Eq.\ \jtildejhat, proving that $\lim_{k\to\infty}T^kJ=\jstar_S$.
   
There remains to
show that $\jstar_S=\jstar $ and that there exists an optimal $S$-regular policy. To this end, we note that by the monotonicity Assumption \assumptionmon, for any policy
$\p=\{\m_0,\m_1,\ldots\}$, we have
$$T_{\m_0}\cdots T_{\m_{k-1}}\bar J\ge T^k\bar J.$$ 
Taking the limit of both sides  as $k\to\infty$, we obtain 
$$J_\p\ge \lim_{k\to\infty}T^k\bar J=\jstar_S,$$
where the equality follows since $T^k J\to \jstar_S$ for all $J\in S$ (as shown earlier), and $\bar J\in S$ [cf.\ Assumption \assumptiontoz(a)].
Thus for all $\p\in \P$, $J_\p\ge \jstar_S=J_{\m},$ 
implying that the policy $\m$ that is optimal within the class of $S$-regular policies is optimal over all policies, and that $\jstar_S=\jstar$. 

\smskip

\pn(c) If  $\m$  is optimal, then $J_\m=\jstar\in S$, so  by Assumption \assumptiontoz(c), $\m$ is
$S$-regular and therefore $T_{\m}J_{\m}=J_{\m}$. Hence, 
$T_{\m}\jstar =T_{\m}J_{\m}=J_{\m}=\jstar =T\jstar.$
Conversely, if
$\jstar =T\jstar  =T_{\m}\jstar $, $\m$ is $S$-regular (cf.\ Lemma \lemmat), so $\jstar=\lim_{k\to\infty}T_{\m}^k\jstar=J_\m$. Therefore,  $\m$ is optimal.
\smskip

\pn(d) If $J\in S$ and $J\le TJ$, by repeatedly applying $T$ to both sides and using the monotonicity of $T$, we obtain
$J\le T^kJ$ for all $k$.
Taking the limit as $k\to\infty$ and using the fact $T^kJ\to \jstar $ [cf.\ part (b)], we obtain $J\le \jstar $. The proof that $J\ge TJ$ implies  $J\ge \jstar$ is similar. 
\smskip

\pn(e) As in  the proof of Prop.\  \propdetpolicyitgen(b), the sequence $\{J_{\m^k}\}$ converges monotonically to a fixed point of $T$, call it $J_\infty$. Since $J_\infty$ lies between $J_{\m^0}\in S$ and $\jstar_S\in S$, it must belong to $S$, by Assumption \assumptiontoz(a). Since the only fixed point of $T$ within $S$ is $\jstar$ [cf.\ part (a)], it follows that $J_\infty=\jstar$.  \qed

Finally let us give a proposition, which provides an approach to verifying part (f) of Assumption \assumptiontoz. The proposition will be used later in this section (cf.\ the proof of Prop.\ 6.4).

\xdef\proppartf{\propn}\propnum\show{myproposition}

\texshopbox{\proposition{\proppartf:}Let $S$ be equal to $R_b(X)$, the subset of $R(X)$ that consists of functions $J$ that are bounded below, i.e., for some $b\in \re$, satisfy $J(x)\ge b$ for all $x\in X$. Let parts (b), (c), and (d) of  Assumption \assumptiontoz\ hold, and assume further that for all scalars $r>0$, we have
$$T\jstar_S-r e\le T(\jstar_S-re),\xdef\partfeq{\lab}\eqnum\show{oneo}$$
where $e$ is the unit function, $e(x)\equiv1$. Then part (f) of Assumption \assumptiontoz\ also holds.}

\proof Let $J\in S$, and let $r>0$ be a scalar such that $\jstar_S-re\le J$ [such a scalar exists since $\jstar_S\in R_b(x)$ by Assumption \assumptiontoz(b)]. Define $J'=\jstar_S-re$, and note that by Lemma \lemmath, $\jstar_S$ is a fixed point of $T$. By using Eq.\ \partfeq, we have
$$J'=\jstar_S-r e=T\jstar_S-r e\le T(\jstar_S-re)=TJ',$$
thus proving part (f) of Assumption \assumptiontoz. \qed

Several examples of applications of Prop.\ \propttffalter\ are given in recent papers of the author. In particular, [Ber15a] considers an application to minimax-type of shortest problems, while [Ber16]  considers an application to SSP problems with multiplicative or exponential cost functions (see also [DeR79], [Pat01], [Ber13], [CaR14]). The paper [Ber15b] considers an infinite-spaces optimal control problem with nonnegative cost per stage, where the objective is to steer a deterministic system towards  a set of termination states. We consider a similar but more general application, where we remove the assumption of nonnegativity for the cost per stage.

\subsubsection{Application to Deterministic Continuous-State Problems}

\pn  Let us consider a deterministic optimal control problem with the system equation
$$x_{k+1}=f(x_k,u_k),\qquad k=0,1,\ldots,\xdef\docsys{\lab}\eqnum\show{oneo}
$$
where $x_k$ and $u_k$ are the state and control at stage $k$, lying in sets $X$ and $U$, respectively, and $f$ is a function mapping $X\times U$ to $X$. The control $u_k$ must be chosen from a constraint set $U(x_k)$. The cost per stage is denoted $g(x,u)$, and is assumed to be a real number.
No restrictions are placed on $X$ and $U$: for example, they may be finite sets as in deterministic shortest path problems, or they may be continuous spaces as in classical problems of control to the origin or some other terminal set.

Because the system is deterministic, given an initial state $x_0$, a policy $\p=\{\m_0,\m_1,\ldots\}$ when applied to the system \docsys, generates a unique sequence of state-control pairs $\big(x_k,\m_k(x_k)\big)$, $k=0,1,\ldots.$ The corresponding cost function is
$$
J_\p(x_0)= \limsup_{N\to\infty} \sum_{k=0}^{N-1}
g\bl(x_k,\mu_k(x_k)\br),\qquad x_0\in X.\eqnum\show{oneo}
$$
We assume that there is a nonempty stopping set $X_0\subset X$, consisting of cost-free and absorbing states in the sense that
$$g(x,u)=0,\qquad x=f(x,u),\qquad \forall\ x\in X_0,\ u\in U(x).\xdef\absorbe{\lab}\eqnum\show{oneo}
$$
Clearly,  for $x\in X_0$, we have $\jstar(x)=0$, as well as  $J_\p(x)=0$ for all policies $\p\in \P$. 
Besides $X_0$, another interesting subset of $X$ is 
$$X_f=\big\{x\in X\mid  \jstar(x)<\infty\big\}.
$$
Ordinarily, in practical applications, the states in $X_f$ are those from which one can reach the stopping set $X_0$, at least asymptotically.

A major class of relevant continuous-state practical problems is control of a dynamic system where the objective is to reach a goal state. Problems of this type are often called {\it planning} problems, and arise frequently in robotics, among others. 
Another major class of practical problems is {\it regulation} problems in control applications, where the objective is to bring and maintain the state within a small region around a desired point. A popular formulation involves a deterministic linear system and a quadratic cost. Variations of this problem may involve a nonquadratic cost function, and state and control constraints.

To formulate a corresponding abstract DP problem, we introduce the mapping $T_\m:R(X)\mapsto R(X)$ by
$$(T_\m J)(x)=g\big(x,\m(x)\big)+J\big(f(x,\m(x))\big),\qquad x\in X,\eqnum\show{oneo}
$$
and the mapping $T:E(X)\mapsto E(X)$ given by
$$(TJ)(x)=\inf_{u\in U(x)}\big\{g(x,u)+J\big(f(x,u)\big)\big\},\qquad x\in X.$$
Here as earlier, we denote by $R(X)$ the set of real-valued functions over $X$, and by $E(X)$ the set of extended real-valued functions over $X$. The initial function $\bar J$ is the zero function [$\bar J(x)\equiv0$]. An important fact is that {\it because the problem is deterministic, $\jstar$ is a fixed point of $T$}.\footnote{$\,$\dag}{\ninepoint  
For any policy $\p=\{\m_0,\m_1,\ldots\}$, using the definition of $J_\p$, we have for all $x$,
$$J_\p(x)=g\big(x,\m_0(x)\big)+J_{\p_1}\big(f(x,\m_0(x))\big),\eqnum\show{oneo}$$
where $\p_1=\{\m_1,\m_2,\ldots\}$. By taking the infimum of the left-hand side over $\p$ and the infimum of the right-hand side over $\p_1$ and then $\m_0$, we obtain $J^*=TJ^*$.}

We say that a policy $\m$ is {\it terminating} if the state sequence $\{x_k\}$ generated starting from any $x\in X_f$ and using $\m$ reaches $X_0$ in finite time, i.e., satisfies $x_{\bar k}\in X_0$ for some index $\bar k$. The set of terminating policies is denoted by ${\cal T}$. Our key assumption is that for $x\in X_f$, the optimal cost $ \jstar(x)$ can be approximated arbitrarily closely by using terminating policies. In particular, we assume the following.

\xdef\assumptiondetoc{\assumptionn}\assumptionnum\show{myproposition}

\texshopbox{\assumption{\assumptiondetoc: (Near-Optimal Termination)}For every pair $(x,\e)$ with $x\in X_f$ and $\e>0$, there exists a terminating policy $\m$ that satisfies $J_\m(x)\le \jstar(x)+\e$.
}
 
This assumption implies in particular that the optimal cost function over terminating policies, 
$$\hat J(x)=\inf_{\m\in{\cal T}}J_\m(x),\qquad x\in X,$$
is equal to $\jstar$. 
Moreover since $\jstar$ is a fixed point of $T$ (because we are dealing with a deterministic problem), it follows that $\hat J$ is a fixed point of $T$, which brings to bear Prop.\ \propregsetcorto.

There are easily verifiable conditions that imply Assumption \assumptiondetoc, some of which are discussed in [Ber15b], where it is assumed in addition that $g\ge0$.  A prominent case is when $X$ and $U$ are finite, so the problem becomes a deterministic shortest path problem.  If all cycles of the state transition graph have positive length, all policies $\p$ that do not terminate from a state $x\in X_f$ must satisfy $J_\p(x)=\infty$, implying that there exists an optimal policy that terminates from all $x\in X_f$. Thus, in this case Assumption \assumptiondetoc\ is naturally satisfied. Another interesting case arises when $g(x,u)=0$ for all $(x,u)$ except if $x\notin X_0$ and $f(x,u)\in X_0$, in which case we have $g(x,u)<0$, i.e., there no cost incurred except for a negative cost (positive reward) upon termination. Then, assuming that $X_0$ can be reached from all states, Assumption \assumptiondetoc\ is satisfied. This is also an example of a deterministic problem where zero length cycles are common. 

When $X$ is the $n$-dimensional Euclidean space $\rn$, a primary case of interest in control system design contexts, it may easily happen that the optimal policies are not terminating from some $x\in X_f$. Instead the optimal state trajectories may approach $X_0$ asymptotically. This is true for example in the classical linear-quadratic optimal control problem, where $X=\rn$, $X_0=\{0\}$, $U=\re^m$, the system is linear of the form $x_{k+1}=Ax_k+Bu_k$, where $A$ and $B$ are given matrices, and the cost is positive semidefinite quadratic. There the optimal policy is linear of the form $\m^*(x)=L x$, where $L$ is some matrix obtained through the steady-state solution of the Riccati equation (see e.g., [Ber05], Section 4.1). Since the optimal closed-loop system is stable and has the form $x_{k+1}=(A+BL)x_k$, the state will typically never reach the termination set $X_0=\{0\}$ in finite time, although it will approach it asymptotically. However, the Assumption \assumptiondetoc\ is satisfied under some natural and easily verifiable controllability and observability conditions (see [Ber15b]). 

Let us denote by $S$ the set of functions
$$S=\big\{J\in E(X)\mid\ J(x)=0,\,\forall\ x\in X_0,\ J(x)\in\re,\,\forall\ x\in X_f,\, J(x)>-\infty,\, \,\forall\ x\in X  \big\}.\xdef\regsets{\lab}\eqnum\show{oneo}
$$
Since $X_0$ consists of cost-free and absorbing states [cf.\  Eq.\ \absorbe], and $\jstar(x)>-\infty$ for all $x\in X$ (by Assumption \assumptiondetoc), the set $S$ contains the cost function $J_\m$ of all policies $\m$, as well as $\jstar$. Moreover it can be seen that every terminating policy is $S$-regular, i.e., ${\cal T}\subset{\cal M}_S$, which implies that
$$\jstar_S=\hat J=\jstar.$$
The reason is that the terminal cost is zero after termination for any terminal cost function $J\in S$, i.e., $(T_\m^kJ)(x)=(T_\m^k\bar J)(x)=J_\m(x)$ for $\m\in{\cal T}$, $x\in X_f$, and $k$ sufficiently large.

The following proposition is a consequence of  Prop.\ \propregsetcorto, the deterministic character of the problem (which guarantees that $\jstar$ is a fixed point of $T$), and Assumption \assumptiondetoc\ (which guarantees that $\jstar_S=\hat J=\jstar$). 

\xdef\propdetoc{\propn}\propnum\show{myproposition}

\texshopbox{\proposition{\propdetoc:}
Let Assumption \assumptiondetoc\ hold. Then:
\nitem{(a)} $\jstar$ is the only  fixed point of $T$ within the set of all $J\in S$ such that $J\ge \jstar$.
\nitem{(b)} We have $T^kJ\to \jstar$ for every
 $J\in S$ such that $J\ge \jstar$.
\nitem{(c)}  If $\m^*$ is terminating and $T_{\m^*}\jstar=T\jstar$, then $\m^*$ is optimal. Conversely, if $\m^*$ is terminating and is optimal, then $T_{\m^*}\jstar=T\jstar$.}

For an example of what may happen in the absence of Assumption \assumptiondetoc, consider the deterministic shortest path Example \exampledetsp\ with $a=0$, $b>0$, and $S=\re$. Here Assumption \assumptiondetoc\ is violated and we have $0=\jstar<\hat J=b$, while the set of fixed points of $T$ is the interval $(-\infty,b]$. 

We will now consider additional assumptions, which guarantee the stronger conclusions of Prop.\ \propttffalter. We first replace the set $S$ of Eq.\ \regsets\ with the following subset of functions that are bounded below:
$$\hat S=\big\{J\in E(X)\mid\ J(x)=0,\,\forall\ x\in X_0,\ J(x)\in\re,\,\forall\ x\in X_f,\, \hbox{$J$ is uniformly bounded below by a scalar} \big\}.
$$
We have the following proposition.

\xdef\propdetoct{\propn}\propnum\show{myproposition}

\texshopbox{\proposition{\propdetoct:}
Let Assumption \assumptiondetoc\ hold, and assume further that:
\nitem{(1)} $\jstar_{\hat S}\in \hat S$.
\nitem{(2)} For each $\hat S$-irregular policy $\m$ and each $J\in \hat S$, there is at least one state $x\in X$ such that  
$\limsup_{k\to\infty}\,(T_\m^kJ)(x)=\infty$.
\nitem{(3)} The control set $U$ is a metric space, and the set
$
\big\{ u\in U(x)\mid g(x,u)+J\big(f(x,u)\big)\le
\l\big\}$
is compact for every
$J\in \hat S$, $x\in X$, and $\l\in \re$.
\pn Then:
\nitem{(a)} The optimal cost function $\jstar$ is the unique fixed point of $T$ within the set $\hat S$.
\nitem{(b)} We have $T^k J\to \jstar$ for all $J\in \hat S$.
\nitem{(c)} A policy  $\m$ is optimal if and only if
$T_\m \jstar =T\jstar$.  Moreover, there exists an optimal $\hat S$-regular policy.
\nitem{(d)} For any $J\in \hat S$, if $J\le TJ$ we have $J\le \jstar$, and if $J\ge TJ$ we have $\hat J\ge \jstar $.
\nitem{(e)} Every sequence $\{\m^k\}$ generated by the PI algorithm starting from an $\hat S$-regular policy $\m^0$ satisfies $J_{\m^k}\downarrow\jstar$.}

\proof The proof consists of showing that all parts of Assumption \assumptiontoz\ are satisfied with $\hat S$ used in place of $S$, so Prop.\ \propttffalter\ applies. Indeed, parts (a) and (e) of this assumption are trivially satisfied, while parts (b)-(d) are the conditions (1)-(3) of the proposition. Then Lemma \lemmath\ is used to assert that $\jstar_{\hat S}$ is a fixed point of $T$. Moreover, Assumption \assumptiontoz(f) is shown using the line of proof of Prop.\ \proppartf. In particular, for any $J\in S$, we let $r>0$ be a scalar such that $\jstar_S-re\le J$ [such a scalar exists since $\jstar_S\in \hat S$ by condition (1)]. Defining $J'=\jstar-re$ where $r>0$ is sufficiently large so that $J'\le J$, we have 
$$J'=\jstar_S-r e=T\jstar_S-r e\le T(\jstar_S-re)=TJ',$$
so Assumption \assumptiontoz(f) holds. Finally the additional assumption needed to apply Prop.\ \propttffalter(e) is clearly satisfied in this deterministic problem. \qed

\vskip-.5pc\vskip-.5pc\vskip-.5pc\vskip-.5pc

\section{Irregular Policies/Finite Cost Case}
\vskip-.5pc

\pn In this section, we consider a perturbation approach to assert that $\jstar_S$ is a fixed point of $T$. This approach applies to problems where some $S$-irregular policies may have finite cost for all states, so Prop.\ \propttffalter\ cannot be used. 
We address this issue by introducing a perturbation that makes the cost of all irregular policies infinite for some states. We can then use  Prop.\ \propdetpolicyitgen\ or Prop.\ \propttffalter\ for the perturbed cost problem, and take the limit as the perturbation vanishes. The idea is that with a perturbation, the cost functions of 
$S$-irregular policies may increase disproportionately relative to the cost functions of the $S$-regular policies, thereby making the problem more amenable to analysis.

In particular, for each $\d\ge0$ and policy $\m$, we consider the mappings $T_{\m,\d}$ and $T_\d$ given by 
$$(T_{\m,\d}J)(x)=H\big(x,\m(x),J\big)+\d,\ \  x\in X,\qquad T_\d J=\inf_{\m\in{\cal M}}T_{\m,\d}J.$$
We define the corresponding cost functions of policies $\p=\{\m_0,\m_1,\ldots\}\in \P$ and $\m\in {\cal M}$, and optimal cost function $\jstar_\d$ by
$$J_{\p,\d}(x)=\limsup_{k\to\infty}\,T_{\m_0,\d}\cdots T_{\m_k,\d}\bar J,\qquad J_{\m,\d}(x)=\limsup_{k\to\infty}\,T_{\m,\d}^k\bar J,\qquad \jstar_\d=\inf_{\p\in \Pi}J_{\p,\d}.$$
We refer to the problem associated with the mappings $T_{\m,\d}$ as the {\it $\d$-perturbed problem\/}. 

The following proposition shows that if the $\d$-perturbed problem is  ``well-behaved" with respect to the $S$-regular policies, then its cost function $\jstar_\d$ can be used to  approximate the optimal cost function $\jstar_S$ over the $S$-regular policies only, and moreover $\jstar_S$ is a fixed point of $T$.

\xdef\proptaffmonf{\propn}\propnum\show{myproposition}

\texshopboxnb{\proposition{\proptaffmonf:} Given a set $S\subset E(X)$, assume that:
\nitem{(1)} For every $\d>0$, we have $\jstar_\d=T_{\d}\jstar_\d$, and there exists an $S$-regular policy $\m^*_\d$ that is optimal for the $\d$-perturbed problem, i.e., $J_{{\m^*_\d},\d}=\jstar_\d$.
\nitem{(2)} For every $S$-regular policy $\m$, we have
$$J_{\m,\d}\le J_\m+w_\m(\d),\qquad \forall\ \d>0,$$
where $w_\m$ is a function such that $\lim_{\d\downarrow0}w_\m(\d)=0$.}\texshopboxnt{\nitem{}
\smskip
\pn Consider $\jstar_S$, the optimal cost function over the $S$-regular policies only: $\jstar_S=\inf_{\m:\, S\hbox{\eightpoint -regular}}J_\m.$
\nitem{(a)} We have $\lim_{\d\downarrow0}\jstar_\d=\jstar_S.$
\nitem{(b)} Assume in addition that $H$ has the property that for every sequence $\{J_m\}\subset S$ with $J_m\downarrow J$, we have
$$\lim_{m\to\infty}H(x,u,J_m)= H(x,u,J),\qquad \forall\ x\in X,\ u\in U(x).\xdef\continuity{\lab}\eqnum\show{oneo}$$
Then $\jstar_S$ is a fixed point of $T$ and the conclusions of Prop.\  \propregsetcorto\ hold.
}

\proof (a) For all $\d>0$, by using conditions (1) and (2), we have for all $S$-regular $\m$,
$$\jstar_S\le J_{\m^*_\d}\le J_{\m^*_\d,\d}= \jstar_\d\le J_{\m,\d}\le J_{\m}+w_{\m}(\d).
$$
By taking the limit as $\d\downarrow0$ and then the infimum over all  $S$-regular $\m$, it follows that
$$\jstar_S\le \lim_{\d\downarrow0}\jstar_\d\le \inf_{\m:\, S\hbox{\eightpoint -regular}}J_\m=\jstar_S.$$
\smskip
\pn (b) From condition (1), for all $\d>0$, we have 
$\jstar_\d= T_\d \jstar_\d\ge T \jstar_\d=TJ_{{\m^*_\d},\d}\ge T\jstar_S,$
and by taking the limit as $\d\downarrow0$ and using part (a), we obtain $\jstar_S\ge T\jstar_S.$
For the reverse inequality, let $\{\d_m\}$ be a sequence with $\d_m\downarrow0$. Using condition (1), we have
$T_{\d_m}\jstar_{\d_m}=\jstar_{\d_m}$, so that for all $m$,
$$H(x,u,\jstar_{\d_m})+\d_m\ge (T_{\d_m}\jstar_{\d_m})(x)=\jstar_{\d_m}(x),\qquad \forall\ x\in X,\ u\in U(x).$$
Taking the limit as $m\to\infty$, and using Eq.\ \continuity\ and the fact $\jstar_{\d_m}\downarrow \jstar_S$ [cf.\ part (a)], we have
$$H(x,u,\jstar_S)\ge \jstar_S(x),\qquad \forall\ x\in X,\ u\in U(x),$$
so that $T\jstar_S\ge \jstar_S$. Thus $\jstar_S$ is a fixed point of $T$, and the assumptions of Prop.\ \propregsetcorto\ are satisfied. \qed

The preceding proposition applies even if 
$\lim_{\d\downarrow0}\jstar_\d(x)>\jstar(x)$
 for some $x\in X$. This is illustrated by the deterministic shortest path Example \exampledetsp, for the zero-cycle case where $a=0$ and $b>0$. Then for $S=\re$, we have $\jstar_S=b>0=\jstar$, while the proposition applies because its assumptions are satisfied. Consistently with the conclusions of the proposition, we have $\jstar_\d=b+\d$, so $\jstar_S= \lim_{\d\downarrow0}\jstar_\d$ and $\jstar_S$ is a fixed point of $T$. We refer to [Ber13] and [BeY16] for further discussion and applications of the approach of this section, and also for a PI algorithm to find $\jstar_S$, which is based on perturbations.

\vskip-1pc

\section{Concluding Remarks}
\vskip-0.5pc

\pn We have provided an analysis of challenging  abstract DP models based on the notion of a regular policy. In particular, we have extended this notion to nonstationary policies, and we have highlighted its connection to an earlier development for stationary policies. We have also streamlined and strengthened the corresponding analysis based on PI-related ideas. The main approach is to start from an interesting set of policy-state pairs satisfying a regularity property, and then characterize the region of convergence of VI. We have shown that this approach can lead to new results in the context of a variety of optimal control problems. In addition to the applications described in this paper, our approach has been applied to minimax and exponential cost shortest path problems [Ber15a], [Ber16].
Our approach may also be applied to other types problems that  involve a termination state and fit the abstract DP framework of this paper, including SSP game problems [PaB99], [Yu11]. These and other related applications are interesting subjects for further research.

\vskip-1pc

\section{References}
\vskip-0.9pc
\def\ref{\vskip1.pt\pn}
\ninepoint

\ref [BeS78]  Bertsekas, D.\ P., and Shreve, S.\ E., 1978.\  Stochastic Optimal
Control:  The Discrete Time Case, Academic Press, N.\ Y.; may be downloaded 
from http://web.mit.edu/dimitrib/www/home.html

\ref [BeT91]  Bertsekas, D.\ P., and Tsitsiklis, J.\ N., 1991.\ ``An Analysis of
Stochastic Shortest Path Problems,"
Math.\ of Operations Research, Vol.\ 16, pp.\ 580-595.

\ref [BeT96]  Bertsekas, D.\ P., and Tsitsiklis, J.\ N., 1996.\ Neuro-Dynamic
Programming, Athena Scientific, Belmont, MA.

\ref[BeY16] Bertsekas, D.\ P., and Yu, H., 2016.\ ``Stochastic Shortest Path Problems Under Weak Conditions,"  Lab.\ for Information and Decision Systems Report LIDS-2909, revision of Jan.\ 2016.

\ref [Ber75] Bertsekas, D.\ P., 1975.\  ``Monotone Mappings in Dynamic Programming," Proc.\ 1975 IEEE Conference on Decision and Control, Houston, TX, pp.\ 20-25. 

\ref [Ber77] Bertsekas, D.\ P., 1977.\  ``Monotone Mappings with Application in
Dynamic Programming," SIAM J.\ on Control and Optimization, Vol.\ 15, pp.\
438-464.

\ref[Ber05] Bertsekas, D.\ P., 2005.\ Dynamic Programming and Optimal Control, Vol.\ I, Athena Scientific, Belmont, MA.

\ref[Ber12] Bertsekas, D.\ P., 2012.\ Dynamic Programming and Optimal Control, Vol.\ II: Approximate Dynamic Programming, Athena Scientific, Belmont, MA.

\ref[Ber13] Bertsekas, D.\ P., 2013.\ Abstract Dynamic Programming, Athena Scientific, Belmont, MA.

\ref[Ber15a] Bertsekas, D.\ P., 2015.\ ``Robust Shortest Path Planning and Semicontractive Dynamic Programming," Lab. for Information and Decision Systems Report LIDS-P-2915, MIT; to appear in Naval Research Logistics.

\ref[Ber15b] Bertsekas, D.\ P., 2015.\ ``Value and Policy Iteration in Deterministic Optimal Control and Adaptive Dynamic Programming," Lab.\ for Information and Decision Systems Report LIDS-P-3174, MIT, Sept. 2015; to appear in IEEE Transactions on Neural Networks and Learning Systems.

\ref[Ber16] Bertsekas, D.\ P., 2016.\ ``Affine Monotonic and Risk-Sensitive Models in Dynamic Programming", Lab. for Information and Decision Systems Report LIDS-3204, MIT. 

\ref [Bla65] Blackwell, D., 1965.\  ``Positive Dynamic Programming," Proc.\ Fifth Berkeley Symposium Math.\ Statistics and Probability, pp.\ 415-418.

\ref[CaR14] Cavus, O., and Ruszczynski, A., 2014.\ ``Risk-Averse Control of Undiscounted Transient Markov Models," SIAM J.\ on Control and Optimization, Vol.\ 52, pp.\ 3935-3966.

\ref[DeR79] Denardo, E.\ V., and Rothblum, U.\ G., 1979.\ ``Optimal Stopping, Exponential Utility, and Linear Programming," Math.\ Programming, Vol.\ 16, pp.\ 228-244.

\ref [Den67] Denardo, E.\ V., 1967.\  ``Contraction Mappings in the Theory Underlying
Dynamic Programming," SIAM Review, Vol.\ 9, pp.\ 165-177.

\ref [Der70] Derman, C., 1970.\ Finite State Markovian Decision Processes,
Academic Press, N.\ Y.

\ref[Fei02] Feinberg, E.\ A., 2002. ``Total Reward Criteria," in E.\ A.\ Feinberg and A.\ Shwartz, (Eds.), Handbook of Markov Decision Processes, Springer, N.\ Y.

\ref[HeL99] Hernandez-Lerma, O., and Lasserre, J.\ B., 1999.\ Further Topics on Discrete-Time Markov Control Processes, Springer, N.\ Y.

\ref[PaB99] Patek, S.\ D., and Bertsekas, D.\ P., 1999.\ ``Stochastic Shortest Path
Games," SIAM J.\ on Control and Opt., Vol.\ 36, pp.\ 804-824.

\ref[Pal67] Pallu de la Barriere, R., 1967.\ Optimal Control Theory, Saunders, Phila; republ.\ Dover, N. Y., 1980.

\ref[Pat01] Patek, S.\ D., 2001.\ ``On Terminating Markov Decision Processes with a Risk Averse Objective Function," Automatica, Vol.\ 37, pp.\ 1379-1386.

\ref [Put94] Puterman, M.\ L., 1994.\  Markov Decision Processes: Discrete Stochastic Dynamic Programming, J.\ Wiley, N.\ Y.

\ref [Str66] Strauch, R., 1966.\  ``Negative Dynamic Programming," Ann.\ Math.\
Statist., Vol.\ 37, pp.\ 871-890.

\ref[Van81]  Van der Wal, J., 1981.\ Stochastic Dynamic Programming, Thesis, The Math.\ Centre, Amsterdam.

\ref [Whi82] Whittle, P., 1982.\  Optimization Over Time, Wiley, N.\ Y., Vol.\
1, 1982, Vol.\ 2, 1983.

\ref[YuB13]  Yu, H.,  and Bertsekas, D.\ P., 2013.\ ``A Mixed Value and Policy Iteration Method for Stochastic Control with Universally Measurable Policies," Lab.\ for Information and Decision Systems Report LIDS-P-2905, MIT; to appear in Math.\ of Operations Research. 

\ref[Yu11] Yu, H., 2011.\ ``Stochastic Shortest Path Games and Q-Learning,"  Lab.\ for Information and Decision Systems Report LIDS-P-2875, MIT.

\ref[Yu15] Yu, H., 2015.\ ``On Convergence of Value Iteration for a Class of Total Cost Markov Decision Processes," SIAM J.\ on Control and Optimization, Vol.\ 53, pp.\ 1982-2016.

\end

%% file: TEXSHOP_macros_new.tex

\def\ignore#1{}
 

\newcount\sectnum
\newcount\subsectnum
\newcount\eqnumber

\global\eqnumber=1\sectnum=0


\def\lab{(\the\sectnum.\the\eqnumber)}



\def\show#1{#1}



\def\smskip{\vskip 5 pt}
\def\medskip{\vskip 10 pt}
\def\bigskip{\vskip 15 pt}
\def\pn{\par\noindent}
\def\br{\break}

\def\bl{\bigl} 
\def\br{\bigr} 
\def\lf{\left}
\def\ri{\right}

\def\tends{\rightarrow}

\def\implies{\Rightarrow}

\def\frac#1#2{{#1\over #2}}

\def\ol#1{\overline{#1}}

\def\a{\alpha}

\def\b{\beta}
\def\l{\lambda}

\def\m{\mu}
\def\p{\pi}

\def\e{\epsilon}

\def\d{\delta}

\def\P{\Pi}

\def\re{\Re}
\def\rn{\Re^n}

\def\tl{\tilde}

\def\old#1{}
\def\leaderfill{\leaders\hbox to 1em{\hss.\hss}\hfill}


\parindent=2pc
\baselineskip=15pt
\vsize=8.7 true in
\voffset=0.125 true in
\parskip=3pt


\def\minprob#1#2#3{$$\eqalign{&\hbox{minimize\ \ }#1\cr &\hbox{subject to\ \
}#2\cr}\ifnum 0=#3{}\else\eqno(#3)\fi$$}        
     
\def\maxprob#1#2#3{$$\eqalign{&\hbox{maximize\ \ }#1\cr &\hbox{subject to\ \
}#2\cr}\ifnum 0=#3{}\else\eqno(#3)\fi$$}        
     
\def\aligntwo#1#2#3#4#5{$$\eqalign{#1&#2\cr #3&#4\cr}
\ifnum 0=#5{}\else\eqno(#5)\fi$$}
\def\alignthree#1#2#3#4#5#6#7{$$\eqalign{#1&#2\cr #3&#4\cr #5&#6\cr}
\ifnum 0=#7{}\else\eqno(#7)\fi$$}


\def\eqnum{\eqno{\hbox{(\the\sectnum.\the\eqnumber)}\global\advance\eqnumber
by1}}

\def\eqnu{\eqno{\hbox{(\the\sectnum.\the\eqnumber)}\global\advance\eqnumber
by1}}

\newcount\examplnumber
\def\examplnum{\global\advance\examplnumber by1}

\newcount\figrnumber
\def\figrnum{\global\advance\figrnumber by1}

\newcount\propnumber
\def\propnum{\global\advance\propnumber by1}

\newcount\defnumber
\def\defnum{\global\advance\defnumber by1}

\newcount\lemmanumber
\def\lemmanum{\global\advance\lemmanumber by1}

\newcount\assumptionnumber
\def\assumptionnum{\global\advance\assumptionnumber by1}

\newcount\conditionnumber
\def\conditionnum{\global\advance\conditionnumber by1}

\def\exampl{\the\sectnum.\the\examplnumber}
\def\figr{\the\sectnum.\the\figrnumber}
\def\propn{\the\sectnum.\the\propnumber}
\def\defn{\the\sectnum.\the\defnumber}
\def\lemman{\the\sectnum.\the\lemmanumber}
\def\assumptionn{\the\sectnum.\the\assumptionnumber}
\def\condn{\the\sectnum.\the\conditionnumber}

\def\section#1{\goodbreak\vskip 3pc plus 6pt minus 3pt\leftskip=-2pc
   \global\advance\sectnum by 1\eqnumber=1
\global\examplnumber=1\figrnumber=1\propnumber=1\defnumber=1\lemmanumber=1\assumptionnumber=1 \conditionnumber =1%
   \line{\hfuzz=1pc{\hbox to 3pc{\bf 
   \vtop{\hfuzz=1pc\hsize=38pc\hyphenpenalty=10000\noindent\uppercase{\the\sectnum.\quad #1}}\hss}}
			\hfill}
			\leftskip=0pc\nobreak\tenf
			\vskip 1pc plus 4pt minus 2pt\noindent\ignorespaces}



\def\sect#1{\noindent\leftskip=-2pc\tenf
   \goodbreak\vskip 1pc plus 4pt minus 2pt
                \global\advance\subsectnum by 1\eqnumber=1
   \line{\hfuzz=1pc{\hbox to 3pc{\bf 
   \vtop{\hfuzz=1pc\hsize=38pc\hyphenpenalty=10000\noindent\uppercase{{\bf #1}}}\hss}}
                        \hfill}
   \leftskip=0pc\nobreak\tenf
                        \vskip 1pc plus 4pt minus 2pt\nobreak\noindent\ignorespaces}

\def\subsection#1{\noindent\leftskip=0pc\tenf
   \goodbreak\vskip 1pc plus 4pt minus 2pt
   \line{\hfuzz=1pc{\hbox to 3pc{\bf 
   \vtop{\hfuzz=1pc\hsize=38pc\hyphenpenalty=10000\noindent{\bf #1}}\hss}}
                        \hfill}
   \leftskip=0pc\nobreak\tenf
                        \vskip 1pc plus 4pt minus 2pt\nobreak\noindent\ignorespaces}
\def\subsubsection#1{\goodbreak\vskip 1pc plus 4pt minus 2pt
   \hfuzz=3pc\leftskip=0pc\noindent\tenit #1 \nobreak\tenf\vskip 6pt plus 1pt
                                minus 1pt\nobreak\ignorespaces\leftskip=0pc}
%

\def\beginexample#1{\noindent\goodbreak\vskip 6pt plus 1pt minus 1pt
\noindent
  \hbox {\bf Example #1\hss}
  \nobreak\vskip 4pt plus 1pt minus 1pt \nobreak\noindent\ninef
  \global\advance
                \leftskip by\parindent\pn}
\def\endexample{\vskip 12pt\tenf\par
  \global\advance\leftskip by -\parindent
  }

\def\beginexercise#1{\noindent\goodbreak\vskip 6pt plus 1pt minus 1pt \noindent\global\normalbaselineskip=12pt
  \hbox {\bf Exercise #1\hss}
  \nobreak\vskip 4pt plus 1pt minus 1pt 
  \nobreak\noindent\ninef\global\advance\leftskip
                        by\parindent\pn}
\def\endexercise{\vskip 12pt\tenf\par
  \global\advance\leftskip by -\parindent
  }

\def\beginsection#1{\noindent\goodbreak\vskip 6pt plus 1pt minus 1pt \noindent\global\normalbaselineskip=12pt
  \hbox {\it #1\hss}
  \vskip 0.1pt plus 1pt minus 1pt \nobreak\noindent\ninef\global\advance
                \leftskip by\parindent\noindent\pn}
\def\endsection{\vskip 12pt\tenf\par
  \global\advance\leftskip by -\parindent
}

%


\def\lemma#1{\smskip\pn{\bf Lemma #1}\quad}

\def\proposition#1{\smskip\pn{\bf Proposition #1}\quad}
\def\proof{\smskip\pn{\bf Proof:}\quad} 
\def\definition#1{\smskip\pn{\bf Definition #1}\quad} 
\def\assumption#1{\smskip\pn{\bf Assumption #1}\quad}

 \def\qed{\quad{\bf
Q.E.D.} \par\bigskip}
\def\ref{\smskip\pn}

\def\chapter#1#2{{\bf \centerline{\helbigbig
{#1}}}\bigskip\bigskip{\bf \centerline{\helbigbig
{#2}}}\bigskip\bigskip} 



\def\longpapertitle#1#2#3{{\bf \centerline{\helbigb
{#1}}}\bigskip{\bf \centerline{\helbigb
{#2}}}\bigskip\bigskip{\centerline{
by}}\bigskip{\bf \centerline{
{#3}}}\bigskip\bigskip} 


\def\nitem#1{\smskip\item{#1}}

\newcount\alphanum
\newcount\romnum

\def\alphaenumerate{\ifcase\alphanum \or (a)\or (b)\or (c)\or (d)\or (e)\or
(f)\or (g)\or (h)\or (i)\or (j)\or (k)\fi}
\def\romenumerate{\ifcase\romnum \or (i)\or (ii)\or (iii)\or (iv)\or (v)\or
(vi)\or (vii)\or (viii)\or (ix)\or (x)\or (xi)\fi}

\def\alist{\begingroup\vskip10pt\alphanum=1
\parskip=2pt\parindent=0pt \leftskip=3pc
\everypar{\llap{{\rm\alphaenumerate\hskip1em}}\advance\alphanum by1}}

\def\nolist{\begingroup\vskip10pt\alphanum=0
\parskip=2pt\parindent=0pt \leftskip=3pc
\everypar{\llap{\global\advance\alphanum by1(\the\alphanum)\hskip1em}}}

\def\romlist{\begingroup\vskip10pt\romnum=1
\parskip=2pt\parindent=0pt \leftskip=5pc
\everypar{\llap{{\rm\romenumerate\hskip1em}}\advance\romnum by1}}



\long\def\fig#1#2#3{\vbox{\vskip1pc\vskip#1
\prevdepth=12pt \baselineskip=12pt
\vskip1pc
\hbox to\hsize{\hfill\vtop{\hsize=25pc\noindent{\eightbf Figure #2\ }
{\eightpoint#3}}\hfill}}}

\long\def\widefig#1#2#3{\vbox{\vskip1pc\vskip#1
\prevdepth=12pt \baselineskip=12pt
\vskip1pc
\hbox to\hsize{\hfill\vtop{\hsize=28pc\noindent{\eightbf Figure #2\ }
{\eightpoint#3}}\hfill}}}

\long\def\table#1#2{\vbox{\vskip0.5pc
\prevdepth=12pt \baselineskip=12pt
\hbox to\hsize{\hfill\vtop{\hsize=25pc\noindent{\eightbf Table #1\ }
{\eightpoint#2}}\hfill}}}

 
\def\rightheadline#1{\headline{\tenrm\hfil #1}}


\long\def\leftfig#1#2{\vbox{\smskip\hsize=220pt
\vtop{{\noindent {\bf #1}}}
\smskip
\noindent
\vbox{{\noindent #2}}
}}

\long\def\rightfig#1#2#3{\vbox{\smskip\vskip#1
\prevdepth=12pt \baselineskip=12pt
\hsize=210pt
\smskip
\vbox{\noindent{\eightbold #2}
\hskip1em{\eightpoint#3}}
}}

\long\def\concept#1#2#3#4#5{\bigskip\hrule
\vbox{\hbox{\leftfig{#1}{#2} \hskip3em
\rightfig{#3}{#4}{#5}} \smskip}
\hrule\bigskip}


\long\def\bconcept#1#2#3#4#5#6#7{
\vbox{
\hbox to \hsize{\vtop{\par #1}}
\concept{#2}{#3}{#4}{#5}{#6}
\hbox to \hsize{\vtop{\par #7}}
\smskip}
}




\def\boxit#1{\vbox{\hrule\hbox{\vrule\kern3pt
                                \vbox{\kern3pt#1\kern3pt}\kern3pt\vrule}\hrule}}
\def\centerboxit#1{$$\vbox{\hrule\hbox{\vrule\kern3pt
                                \vbox{\kern3pt#1\kern3pt}\kern3pt\vrule}\hrule}$$}

\long\def\boxtext#1#2{$$\boxit{\vbox{\hsize #1\noindent\strut #2\strut}}$$}

%
%
%

\def\picture #1 by #2 (#3){
  \vbox to #2{
    \hrule width #1 height 0pt depth 0pt
    \vfill
    \special{picture #3} 
    }
  }

\def\scaledpicture #1 by #2 (#3 scaled #4){{
  \dimen0=#1 \dimen1=#2
  \divide\dimen0 by 1000 \multiply\dimen0 by #4
  \divide\dimen1 by 1000 \multiply\dimen1 by #4
  \picture \dimen0 by \dimen1 (#3 scaled #4)}
  }

%
%

\long\def\captfig#1#2#3#4#5{\vbox{\vskip1pc
\hbox to\hsize{\hfill{\picture #1 by #2 (#3)}\hfill}
\prevdepth=9pt \baselineskip=9pt
\vskip1pc
\hbox to\hsize{\hfill\vtop{\hsize=24pc\noindent{\eightbold Figure #4}
\hskip1em{\eightpoint#5}}\hfill}}}

%
%
%

\def\illustration #1 by #2 (#3){
  \vskip#2\hskip#1\special{illustration #3} 
    }

\def\scaledillustration #1 by #2 (#3 scaled #4){{
  \dimen0=#1 \dimen1=#2
  \divide\dimen0 by 1000 \multiply\dimen0 by #4
  \divide\dimen1 by 1000 \multiply\dimen1 by #4
  \illustration \dimen0 by \dimen1 (#3 scaled #4)}
  }


\newbox\graybox
\newdimen\xgrayspace
\newdimen\ygrayspace
%
%
%
%
%
%
%
%
%

\def\Textshade#1#2#3#4#5#6{%
    \xgrayspace=#4pt%
    \ygrayspace=#4pt%
    \def\grayshade{#3}%
    \def\linewidth{#5}%
    \def\theradius{#6}%
    \setbox\graybox=\hbox{\surroundboxa{#2}}%
    \hbox{%
    \hbox to 0pt{%
    \PScommands
}%
    \box\graybox}}%
%
%
\long%

\long%
\def\Parashade#1#2#3#4#5#6#7{%
    \xgrayspace=#4pt%
    \ygrayspace=#4pt%
    \def\grayshade{#3}%
    \def\linewidth{#5}%
    \def\theradius{#6}%
    \def\thevskip{#7pt}%
    \setbox\graybox=\hbox{\surroundboxb{#2}}%
    \vskip\thevskip%
    \hbox{%
    \hbox to 0pt{%
    \PScommands
}%
     \box\graybox}%
     \vskip\thevskip%
}%
%
%
%
\long\def\surroundboxa#1{\leavevmode\hbox{\vtop{%
\vbox{\kern\ygrayspace%
\hbox{\kern\xgrayspace#1%
      \kern\xgrayspace}}\kern\ygrayspace}}}
%
%
\long\def\surroundboxb#1{\leavevmode\hbox{\vtop{%
\vbox{\kern\ygrayspace%
\hbox{\kern\xgrayspace\vbox{\advance\hsize-2\xgrayspace#1}%
      \kern\xgrayspace}}\kern\ygrayspace}}}
%
%
%
\long\def\PScommands{%
\special{rawpostscript
/sharpbox{%
           newpath
           xmin ymin moveto
           xmin ymax lineto
           xmax ymax lineto
           xmax ymin lineto
           xmin ymin lineto
           closepath 
          }bind def
}%
\special{rawpostscript
/sharpboxnb{%
           newpath
           xmin ymin moveto
           xmin ymax lineto
           xmax ymax lineto
           xmax ymin lineto
          }bind def
}%
\special{rawpostscript
/sharpboxnt{%
           newpath
           xmin ymax moveto
           xmin ymin lineto
           xmax ymin lineto
           xmax ymax lineto
          }bind def
}%
\special{rawpostscript
/roundbox{%
           newpath
           xmin radius add ymin moveto
           xmax ymin xmax ymax radius arcto
           xmax ymax xmin ymax radius arcto
           xmin ymax xmin ymin radius arcto
           xmin ymin xmax ymin radius arcto 16 {pop} repeat
           closepath
          }bind def
}%
\special{rawpostscript
/sharpcorners{%
               sharpbox gsave grayshade setgray fill grestore 
               linewidth setlinewidth stroke
              }bind def
}%
\special{rawpostscript
/sharpcornersnt{%
               sharpboxnt gsave grayshade setgray fill grestore 
               linewidth setlinewidth stroke
              }bind def
}%
\special{rawpostscript
/sharpcornersnb{%
               sharpboxnb gsave grayshade setgray fill grestore 
               linewidth setlinewidth stroke
              }bind def
}%
\special{rawpostscript
/roundcorners{%
               roundbox gsave grayshade setgray fill grestore 
               linewidth setlinewidth stroke
              }bind def
}%
\special{rawpostscript
/plainbox{%
           sharpbox grayshade setgray fill 
          }bind def
}%
%
\special{rawpostscript
/roundnoframe{%
               roundbox grayshade setgray fill 
              }bind def
}%
\special{rawpostscript
/sharpnoframe{%
               sharpbox grayshade setgray fill 
              }bind def
}%
}%
%
%

\def\pshade#1{\Parashade{sharpcorners}{#1}{0.95}{10}{0.5}{10}{10}}


\def\boxit#1{\vbox{\hrule\hbox{\vrule\kern3pt
                                \vbox{\kern3pt#1\kern3pt}\kern3pt\vrule}\hrule}}

\def\boxitnb#1{\vbox{\hrule\hbox{\vrule\kern3pt
                                \vbox{\kern3pt#1\kern3pt}\kern3pt\vrule}}}

\def\boxitnt#1{\vbox{\hbox{\vrule\kern3pt
                                \vbox{\kern3pt#1\kern3pt}\kern3pt\vrule}\hrule}}

\long\def\boxtext#1#2{$$\boxit{\vbox{\hsize #1\noindent\strut #2\strut}}$$}
\long\def\boxtextnb#1#2{$$\boxitnb{\vbox{\hsize #1\noindent\strut #2\strut}}$$}
\long\def\boxtextnt#1#2{$$\boxitnt{\vbox{\hsize #1\noindent\strut #2\strut}}$$}

\def\texshopbox#1{\boxtext{462pt}{\vskip-1.5pc\pshade{\vskip-1.0pc#1\vskip-2.0pc}}}
\def\texshopboxnt#1{\boxtextnt{462pt}{\vskip-1.5pc\pshade{\vskip-1.0pc#1\vskip-2.0pc}}}
\def\texshopboxnb#1{\boxtextnb{462pt}{\vskip-1.5pc\pshade{\vskip-1.0pc#1\vskip-2.0pc}}}


%
%
%
%
%
%
%
%
\font\helbigbig=cmr10 scaled 2500%
\font\helbigb=cmbx10 scaled 1500%
\font\eightbold=cmbx8%

\def\tenf{\hel}%
\def\tenit{\heli}%
\def\ninef{\ninehel}%
\def\nineit{\nineheli}%
%
%


\font\tenrm=cmr10%
\font\teni=cmmi10%
\font\tensy=cmsy10%
\font\tenbf=cmbx10%
\font\tentt=cmtt10%
\font\tenit=cmti10%
\font\tensl=cmsl10%

\def\tenpoint{\def\rm{\fam0\tenrm}%
\textfont0=\tenrm%
\textfont1=\teni%
\textfont2=\tensy%
\textfont\itfam=\tenit%
\textfont\slfam=\tensl%
\textfont\ttfam=\tentt%
\textfont\bffam=\tenbf%
\scriptfont0=\sevenrm%
\scriptfont1=\seveni%
\scriptfont2=\sevensy%
\scriptscriptfont0=\sixrm%
\scriptscriptfont1=\sixi%
\scriptscriptfont2=\sixsy%
\def\it{\fam\itfam\tenit}%
\def\tt{\fam\ttfam\tentt}%
\def\sl{\fam\slfam\tensl}%
\scriptfont\bffam=\sevenbf%
\scriptscriptfont\bffam=\sixbf%
\def\bf{\fam\bffam\tenbf}%
\normalbaselineskip=18pt%
\normalbaselines\rm}%

\font\ninerm=cmr9%
\font\ninebf=cmbx9%
\font\nineit=cmti9%
\font\ninesy=cmsy9%
\font\ninei=cmmi9%
\font\ninett=cmtt9%
\font\ninesl=cmsl9%

\def\ninepoint{\def\rm{\fam0\ninerm}%
\textfont0=\ninerm%
\textfont1=\ninei%
\textfont2=\ninesy%
\textfont\itfam=\nineit%
\textfont\slfam=\ninesl%
\textfont\ttfam=\ninett%
\textfont\bffam=\ninebf%
\scriptfont0=\sixrm%
\scriptfont1=\sixi%
\scriptfont2=\sixsy%
\def\it{\fam\itfam\nineit}%
\def\tt{\fam\ttfam\ninett}%
\def\sl{\fam\slfam\ninesl}%
\scriptfont\bffam=\sixbf%
\scriptscriptfont\bffam=\fivebf%
\def\bf{\fam\bffam\ninebf}%
\normalbaselineskip=16pt%
\normalbaselines\rm}%

\font\eightrm=cmr8%
\font\eighti=cmmi8%
\font\eightsy=cmsy8%
\font\eightbf=cmbx8%
\font\eighttt=cmtt8%
\font\eightit=cmti8%
\font\eightsl=cmsl8%

\def\eightpoint{\def\rm{\fam0\eightrm}%
\textfont0=\eightrm%
\textfont1=\eighti%
\textfont2=\eightsy%
\textfont\itfam=\eightit%
\textfont\slfam=\eightsl%
\textfont\ttfam=\eighttt%
\textfont\bffam=\eightbf%
\scriptfont0=\sixrm%
\scriptfont1=\sixi%
\scriptfont2=\sixsy%
\scriptscriptfont0=\fiverm%
\scriptscriptfont1=\fivei%
\scriptscriptfont2=\fivesy%
\def\it{\fam\itfam\eightit}%
\def\tt{\fam\ttfam\eighttt}%
\def\sl{\fam\slfam\eightsl}%
\scriptscriptfont\bffam=\fivebf%
\def\bf{\fam\bffam\eightbf}%
\normalbaselineskip=14pt%
\normalbaselines\rm}%

\font\sevenrm=cmr7%
\font\seveni=cmmi7%
\font\sevensy=cmsy7%
\font\sevenbf=cmbx7%

\def\sevenpoint{%
   \def\rm{\sevenrm}\def\bf{\sevenbf}%
   \def\smc{\sevensmc}\baselineskip=12pt\rm}%

\font\sixrm=cmr6%
\font\sixi=cmmi6%
\font\sixsy=cmsy6%
\font\sixbf=cmbx6%

\fontdimen13\tensy=2.6pt%
\fontdimen14\tensy=2.6pt%
\fontdimen15\tensy=2.6pt%
\fontdimen16\tensy=1.2pt%
\fontdimen17\tensy=1.2pt%
\fontdimen18\tensy=1.2pt%

\def\tenf{\tenpoint}%
\def\ninef{\ninepoint}%
%


%% file: TEXSHOP_small_baseline.tex

\def\tenpoint{\def\rm{\fam0\tenrm}%
\textfont0=\tenrm%
\textfont1=\teni%
\textfont2=\tensy%
\textfont\itfam=\tenit%
\textfont\slfam=\tensl%
\textfont\ttfam=\tentt%
\textfont\bffam=\tenbf%
\scriptfont0=\sevenrm%
\scriptfont1=\seveni%
\scriptfont2=\sevensy%
\scriptscriptfont0=\sixrm%
\scriptscriptfont1=\sixi%
\scriptscriptfont2=\sixsy%
\def\it{\fam\itfam\tenit}%
\def\tt{\fam\ttfam\tentt}%
\def\sl{\fam\slfam\tensl}%
\scriptfont\bffam=\sevenbf%
\scriptscriptfont\bffam=\sixbf%
\def\bf{\fam\bffam\tenbf}%
\normalbaselineskip=14pt%
\normalbaselines\rm}%

\def\ninepoint{\def\rm{\fam0\ninerm}%
\textfont0=\ninerm%
\textfont1=\ninei%
\textfont2=\ninesy%
\textfont\itfam=\nineit%
\textfont\slfam=\ninesl%
\textfont\ttfam=\ninett%
\textfont\bffam=\ninebf%
\scriptfont0=\sixrm%
\scriptfont1=\sixi%
\scriptfont2=\sixsy%
\def\it{\fam\itfam\nineit}%
\def\tt{\fam\ttfam\ninett}%
\def\sl{\fam\slfam\ninesl}%
\scriptfont\bffam=\sixbf%
\scriptscriptfont\bffam=\fivebf%
\def\bf{\fam\bffam\ninebf}%
\normalbaselineskip=13pt%
\normalbaselines\rm}%

\def\eightpoint{\def\rm{\fam0\eightrm}%
\textfont0=\eightrm%
\textfont1=\eighti%
\textfont2=\eightsy%
\textfont\itfam=\eightit%
\textfont\slfam=\eightsl%
\textfont\ttfam=\eighttt%
\textfont\bffam=\eightbf%
\scriptfont0=\sixrm%
\scriptfont1=\sixi%
\scriptfont2=\sixsy%
\scriptscriptfont0=\fiverm%
\scriptscriptfont1=\fivei%
\scriptscriptfont2=\fivesy%
\def\it{\fam\itfam\eightit}%
\def\tt{\fam\ttfam\eighttt}%
\def\sl{\fam\slfam\eightsl}%
\scriptscriptfont\bffam=\fivebf%
\def\bf{\fam\bffam\eightbf}%
\normalbaselineskip=12pt%
\normalbaselines\rm}%

\def\sevenpoint{%
   \def\rm{\sevenrm}\def\bf{\sevenbf}%
   \def\smc{\sevensmc}\baselineskip=10pt\rm}%